\documentclass[12pt]{amsart}
\usepackage{amssymb}
\usepackage{geometry}
\usepackage[vcentermath]{youngtab}
\usepackage[all]{xy}
\input prepictex
\input pictexwd
\input postpictex

\newtheorem{theorem}{Theorem}[section]
\newtheorem{lemma}[theorem]{Lemma}
\newtheorem{prop}[theorem]{Proposition}
\newtheorem{corollary}[theorem]{Corollary}
\theoremstyle{definition}

\newtheorem{example}[theorem]{Example}
\newtheorem{remark}[theorem]{Remark}

\newcommand\yone{\raise1.5pt\hbox{$v_1$}}
\newcommand\ytwo{\raise1.5pt\hbox{$v_2$}}
\newcommand\ythree{\raise1.5pt\hbox{$v_3$}}
\newcommand\yfour{\raise1.5pt\hbox{$v_4$}}
\newcommand\yfive{\raise1.5pt\hbox{$v_5$}}
\newcommand\ysix{\raise1.5pt\hbox{$v_6$}}
\newcommand\yseven{\raise1.5pt\hbox{$v_7$}}
\newcommand\yeight{\raise1.5pt\hbox{$v_8$}}
\newcommand\ynine{\raise1.5pt\hbox{$v_9$}}

\title
[Rook placements and Jordan forms]
{Rook placements and Jordan forms of upper-triangular nilpotent matrices}
\author{Martha Yip
}
\email{martha.yip@uky.edu}
\address{Department of Mathematics, University of Kentucky, Lexington KY 40506}

\begin{document}
\begin{abstract}
The set of $n$ by $n$ upper-triangular nilpotent matrices with entries in a finite field $\mathbb{F}_q$ has Jordan canonical forms indexed by partitions $\lambda \vdash n$.  We present a combinatorial formula for computing the number $F_\lambda(q)$ of matrices of Jordan type $\lambda$ as a weighted sum over standard Young tableaux. We also study a connection between these matrices and non-attacking rook placements, which leads to a refinement of the formula for $F_\lambda(q)$.
\end{abstract}
\maketitle 
\parskip=8pt

\section{Introduction}\label{sec:intro}


In the beautiful paper {\em Variations on the Triangular Theme}~\cite{kirillov95},
Kirillov studied various structures on the set of triangular matrices.
Let $G=G_n(\mathbb{F}_q)$ denote the group of $n$ by $n$ invertible upper-triangular matrices over the field $\mathbb{F}_q$ having $q$ elements, and let $\mathfrak{g}=\mathfrak{g}_n(\mathbb{F}_q)=\mathrm{Lie}(G_n(\mathbb{F}_q))$ denote the corresponding Lie algebra of $n$ by $n$ upper-triangular nilpotent matrices over $\mathbb{F}_q$.  
The problem of determining the set $\mathcal{O}_n(\mathbb{F}_q)$ of adjoint $G$-orbits in $\mathfrak{g}$ remains challenging, and a more tractable task is to study a decomposition of $\mathcal{O}_n(\mathbb{F}_q)$ via the Jordan canonical form.
Let $\lambda \vdash n$ be a partition of $n$ with $r$ positive parts $\lambda_1\geq \lambda_2 \geq \cdots \geq \lambda_r>0$, and let
$$J_\lambda = J_{\lambda_1}\oplus J_{\lambda_2}\oplus \cdots \oplus J_{\lambda_r},
\qquad\hbox{where}\qquad
J_i=\begin{bmatrix}
0&1&0 & \cdots & 0& 0 \\ 
0&0&1 & \cdots & 0 &0\\ 
0&0&0 & \cdots & 0 &0\\ 
\vdots& \vdots &\vdots & \ddots & \vdots &\vdots \\
0&0&0 & \cdots &0& 1 \\  
0&0&0 & \cdots & 0&0 \\ 
\end{bmatrix}_{i\times i}
$$
is the $i$ by $i$ elementary Jordan matrix with all eigenvalues equal to zero. 
If $X\in \mathfrak{g}_n(\mathbb{F}_q)$ is similar to $J_\lambda$ under $GL_n(\mathbb{F}_q)$, then $X$ is said to have {\em Jordan type $\lambda$}.
Each conjugacy class contains a unique Jordan matrix $J_\lambda$, so these classes are indexed by the partitions of $n$.
Evidently, the Jordan type of $X$ depends only on its adjoint $G$-orbit.

Let $\mathfrak{g}_{n,\lambda}(\mathbb{F}_q)\subseteq \mathfrak{g}_{n}(\mathbb{F}_q)$ be the set of upper-triangular nilpotent matrices of fixed Jordan type $\lambda$, and let
\begin{equation}\label{eq.Flambda}
F_\lambda(q) = \left|\,\mathfrak{g}_{n,\lambda}(\mathbb{F}_q)\,\right|.
\end{equation}

Springer showed that $\mathfrak{g}_{n,\lambda}(\mathbb{F}_q)$ is an algebraic manifold 
with $f^\lambda$ irreducible components, where $f^\lambda$ is the number of standard Young tableaux of shape $\lambda$, and each of which has dimension ${n\choose 2}-n_\lambda$, where $n_\lambda$ is an integer defined in Equation~\ref{eq.nlambda}.   These quantities appear in the study of $F_\lambda(q)$.

In Section 2, we show that the numbers $F_\lambda(q)$ satisfy a simple recurrence equation, and that they are polynomials in $q$ with integer coefficients.  
As a consequence of the recurrence equation in Theorem~\ref{thm.JFrec}, it follows that the coefficient of the highest degree term in $F_\lambda(q)$ is $f^\lambda$, and
$\deg F_\lambda(q) = {n\choose 2}-n_\lambda$.
Equation~\eqref{eqn.syt} is a combinatorial formula for $F_\lambda(q)$ as a sum over standard Young tableaux of shape $\lambda$ that can be derived from the recurrence equation.

The cases $F_{(1^n)}(q) = 1$ and $F_{(n)}(q) = (q-1)^{n-1} q^{n-1\choose 2}$ are readily computed, since the matrix in $\mathfrak{g}_n(\mathbb{F}_q)$ of Jordan type $(1^n)$ is the matrix of zero rank, and the matrices in $\mathfrak{g}_n(\mathbb{F}_q)$ of Jordan type $(n)$ are the matrices of rank equal to $n-1$. 
Section 2 concludes with explicit formulas for $F_\lambda(q)$ in several other special cases of $\lambda$, including hook shapes, two-rowed partitions and two-columned partitions.

In Section 3, we explore a connection of $F_\lambda(q)$ with rook placements.
In their study of a formula of Frobenius, Garsia and Remmel~\cite{garsiaRemmel86} introduced the {\em $q$-rook polynomial}
\begin{displaymath}
R_{B,k}(q)= \sum_{c\in \mathcal{C}(B,k)} q^{\mathrm{inv}(c)},
\end{displaymath}
which is a sum over the set $\mathcal{C}(B,k)$ of non-attacking placements of $k$ rooks on a Ferrers board $B$, and $\mathrm{inv}(c)$, defined in Equation~\eqref{eq.inv}, is the number of inversions of $c$.  In the case when $B=B_n$ is the staircase-shaped board, Garsia and Remmel showed that
$R_{B_n,k}(q) = S_{n,n-k}(q)$
is a {\em $q$-Stirling number of the second kind}.  These numbers are defined by the recurrence equation
\begin{displaymath}
S_{n,k}(q) = q^{k-1} S_{n-1,k-1}(q) +[k]_qS_{n-1,k}(q) \quad\hbox{for}\quad
0\leq k\leq n,
\end{displaymath}
with initial conditions $S_{0,0}(q)=1$, and $S_{n,k}(q) = 0$ for $k<0$ or $k>n$.

It was shown by Solomon~\cite{solomon90} that 
non-attacking placements of $k$ rooks on rectangular $m\times n$ boards are naturally associated to $m$ by $n$ matrices with rank $k$ over $\mathbb{F}_q$.  
By identifying a Ferrers board $B$ inside an $n$ by $n$ grid with the entries of an $n$ by $n$ matrix, Haglund~\cite{haglund98} generalized Solomon's result to the case of non-attacking placements of $k$ rooks on Ferrers boards, and obtained a formula for the number of $n$ by $n$ matrices with rank $k$ whose support is contained in the Ferrers board region.
A special case of Haglund's formula shows that the number of $n$ by $n$ nilpotent upper-triangular matrices of rank $k$ is
\begin{equation}\label{eq.Pkq1}
P_{B_n,k}(q)
= (q-1)^{k}q^{{n\choose 2}-k}R_{B_n,k}(q^{-1}).
\end{equation}

Now, a matrix in $\mathfrak{g}_{n,\lambda}(\mathbb{F}_q)$ has rank $n-\ell(\lambda)$, 
where $\ell(\lambda)$ is the number of parts of $\lambda$, so the number of matrices in $\mathfrak{g}_n(\mathbb{F}_q)$ with rank $k$ is
\begin{equation}\label{eq.Pkq2}
P_{B_n,k}(q)=\sum_{\lambda\vdash n:\ \ell(\lambda)=n-k} F_\lambda(q).
\end{equation}
Given Equations~\ref{eq.Pkq1} and~\ref{eq.Pkq2}, it would be natural to ask whether it is possible to partition the placements $\mathcal{C}(B_n,k)$ into disjoint subsets so that the sum over each subset of placements gives $F_\lambda(q)$.  A central goal of this paper is to study the connection between upper-triangular nilpotent matrices over $\mathbb{F}_q$ and non-attacking rook placements on the staircase-shaped board $B_n$.  Theorem~\ref{thm.Phi} shows that there is a weight-preserving bijection $\Phi$ between rook placements on $B_n$ and paths in a graph $\mathcal{Z}$ (see Figure~\ref{fig.Z}), which is a multi-edged version of Young's lattice.
As a result, we obtain Corollary~\ref{cor.main}, which gives a formula for $F_\lambda(q)$ as a sum over certain rook placements that can be viewed as a generalization of Haglund's formula in Equation~\eqref{eq.Pkq1}.

There is a classically known bijection between rook placements in $\mathcal{C}(B_n,k)$ and set partitions of $[n]$ with $n-k$ parts, so it is logical to next study the connection between $F_\lambda(q)$ and set partitions.  We do this in Section 4.  Theorem~\ref{thm.Psi} describes the construction of a new (weight-preserving) bijection $\Psi$ between rook placements and set partitions.  These bijections allow us to refine Equation~\eqref{eqn.syt} to a sum over set partitions (or rook placements). We also discuss the significance of the polynomials $F_C(q)$ indexed by rook placements in a special case.

{\bf Acknowledgements.} I would like to thank Jim Haglund and Alexandre Kirillov for their invaluable guidance, and Yue Cai and Alejandro Morales for enlightening conversations.

\section{Formulas for $F_\lambda(q)$}
\label{sec.RecEq}
The recurrence equation for $F_\lambda(q)$ in Theorem~\ref{thm.JFrec} can be found in~\cite[Division Theorem]{borodin95}, where Borodin considers the matrices as particles of a certain mass and studies the asymptotic behaviour of the formula.  A preliminary version of the idea first appeared in~\cite{kirillov95}. In this section, we give an elementary proof of the formula, and investigate some of the combinatorial properties of $F_\lambda(q)$.

\subsection{The recurrence equation for $F_\lambda(q)$}
A {\em partition} $\lambda$ of a nonnegative integer $n$, denoted by $\lambda \vdash n$, is a non-increasing sequence of nonnegative integers $\lambda_1\geq \lambda_2 \geq \cdots \geq \lambda_n \geq 0$ with $|\lambda | = \sum_{i=1}^n \lambda_i=n$.  If $\lambda$ has $r$ positive parts, write $\ell(\lambda)=r$.  A
partition $\lambda$ can be represented by its Ferrers diagram in the English notation, which is an array of $\lambda_i$ boxes in the $i$th row, with the boxes justified upwards and to the left.  Let $\lambda_j'$ denote the size of the $j$th column of $\lambda$.

Young's lattice $\mathcal{Y}$ is the lattice of partitions ordered by the inclusion of their Ferrers diagrams; that is, $\mu \leq \lambda$ if and only if $\mu_i\leq \lambda_i$ for every $i$.  In particular, $\mu$ is covered by $\lambda$ in the Hasse diagram of $\mathcal{Y}$ and we write $\mu \lessdot \lambda$
if the Ferrers diagram of $\lambda$ can be obtained by adding a box to the Ferrers diagram of $\mu$.  See Figure 1.

\begin{example}\Yboxdim10pt
The partition
\begin{displaymath}
\lambda = (4,2,2,1) \vdash 9 
\quad \hbox{ has diagram }\quad  
\yng(4,2,2,1)
\end{displaymath}
and columns $\lambda_1'=4, \lambda_2'=3, \lambda_3'=1, \lambda_4'=1$.
\end{example}

\begin{lemma}\label{lem.rank}
Let $\lambda\vdash n$ be a partition whose Ferrers diagram has $r$ rows and $c$ columns.  The Jordan matrix $J_\lambda$ satisfies
$$\mathrm{rank}\left(J_\lambda^k\right) 
= \begin{cases}
\lambda_{k+1}'+\cdots +\lambda_c', 
	&\hbox{if } 0 \leq k < c,\\ 
0, 
	& \hbox{if } k\geq c.
\end{cases}$$
\end{lemma}
\proof The $i$ by $i$ elementary Jordan matrix $J_i$ has
$\mathrm{rank}\left(J_i^k\right) = i-k$ if $0 \leq k \leq i$, and its rank is zero otherwise, so the Jordan matrix $J_\lambda=J_{\lambda_1}\oplus\cdots\oplus J_{\lambda_r}$ has
$$
\mathrm{rank}\left(J_\lambda^k\right)
= \sum_{i=1}^r \mathrm{rank}\left( J_{\lambda_i}^k \right)
= \sum_{i: \lambda_i\geq k}\mathrm{rank}\left( J_{\lambda_i}^k \right)
= \sum_{j=k+1}^c \lambda_j',
$$
for $0\leq k <c$, which is the number of boxes in the last $c-k$ columns of $\lambda$.
\qed

\begin{remark}\label{rem.reduce}
Matrices which are similar have the same rank, so if $X\in \mathfrak{g}_{n,\lambda}(\mathbb{F}_q)$, then
$\mathrm{rank}\left(X^k\right) = \mathrm{rank}\left(J_\lambda^k\right)$
for all $k\geq0$.
Conversely, let $\lambda, \nu \vdash n$.  It follows from Lemma~\ref{lem.rank} that $\mathrm{rank}\left(J_\lambda^k\right) = \mathrm{rank}\left(J_\nu^k\right)$ for all $k\geq0 $ if and only if $\lambda = \nu$.  Thus if $X\in \mathfrak{g}_n(\mathbb{F}_q)$ is a matrix such that $\mathrm{rank}\left(X^k\right) = \mathrm{rank}\left(J_\lambda^k\right)$ for all $k\geq0$, then $X$ is similar to $J_\lambda$.
\end{remark}

\begin{example}
If a matrix $X\in \mathfrak{g}_n(\mathbb{F}_q)$ has Jordan type $\lambda=(4,2,2,1)$, then $\mathrm{rank}(X) = 5$, $\mathrm{rank}(X^2)=2$, $\mathrm{rank}(X^3)=1$, and $\mathrm{rank}(X^4)=0$.
\end{example}


If $X\in \mathfrak{g}_{n}(\mathbb{F}_q)$ is a matrix of the form
$$X = 
\begin{bmatrix} 
J_\mu & \mathbf{v}\\
\mathbf{0} &0
\end{bmatrix},
$$
where $\mu\vdash n-1$, and $\mathbf{v} = [v_1, \ldots, v_{n-1}]^T \in \mathbb{F}_q^{n-1}$, then the first order leading principal submatrix of $X^k$ is $J_\mu^k$, and for $1\leq k \leq n$, we define column vectors $\mathbf{v}^k = [v_1^k,\ldots, v_{n-1}^k]^T \in \mathbb{F}_q^{n-1}$ by
$$X^k = 
\begin{bmatrix} 
J_\mu^k & \mathbf{v}^k\\
\mathbf{0} &0
\end{bmatrix}.
$$
For $i\geq1$, let $\alpha_i = \mu_1 + \cdots + \mu_i$ be the sum of the first $i$ parts of $\mu$.
The $(i,j)$th entry of $J_\mu^k$ is nonzero if and only if $j=i+k$, and $i,i+k \leq \alpha_b$ for all $b\geq1$.
It follows from this that
\begin{equation}\label{eqn.visualizev}
v_i^k = \begin{cases}
v_{i+k-1}, 
	&\hbox{if $i, i+k-1\leq \alpha_b$ for all $b\geq 1$},\\
0,
	&\hbox{otherwise}.
\end{cases}
\end{equation}

There is a simple way to visualize the vectors $\mathbf{v}^k$, which we illustrate with an example.
\begin{example} \label{eg.visualize}
Let $\mu = (4,2,1,1)$, so that $\alpha_1 = 4, \alpha_2=6, \alpha_3=7$, and $\alpha_4=8$.  Let
$$
X=\begin{bmatrix}
0&1&0&0&&&&&v_1\\
0&0&1&0&&&&&v_2\\
0&0&0&1&&&&&v_3\\
0&0&0&0&&&&&v_4\\ \hline
&&&&0&1&&&v_5\\
&&&&0&0&&&v_6\\ \hline
&&&&&&0&&v_7\\ \hline
&&&&&&&0&v_8\\ \hline
&&&&&&&&0
\end{bmatrix}
\quad\hbox{so that}\quad
X^2=\begin{bmatrix}
0&0&1&0&&&&&v_2\\
0&0&0&1&&&&&v_3\\
0&0&0&0&&&&&v_4\\
0&0&0&0&&&&&0\\ \hline
&&&&0&0&&&v_6\\
&&&&0&0&&&0\\ \hline
&&&&&&0&&0\\ \hline
&&&&&&&0&0\\ \hline
&&&&&&&&0
\end{bmatrix}.
$$
We may visualize the vectors $\mathbf{v}$ and $\mathbf{v}^2$ as fillings of the Ferrers diagram for $\mu$:
\Yboxdim16pt
$$\mathbf{v} 
=\young(\yfour\ythree\ytwo\yone,\ysix\yfive,\yseven,\yeight)
\qquad\hbox{and}\qquad
\mathbf{v}^2 
= \young(0\yfour\ythree\ytwo,0\ysix,0,0).$$
This way, a basis of $\ker X^k$ is the set of vectors filling the first $k$ columns of the diagram.
\end{example}

\begin{lemma}\label{lem.1}
If $X\in \mathfrak{g}_{n,\lambda}(\mathbb{F}_q)$ and its
first order leading principal submatrix $Y\in \mathfrak{g}_{n-1,\mu}(\mathbb{F}_q)$, then $\lambda\gtrdot\mu$.
\end{lemma}
\proof
We first consider the case $Y=J_\mu$.  
If $\mu$ has $s$ parts, let $\alpha_i = \mu_1+\cdots+\mu_i$ for $1\leq i\leq s$.
Then 
\begin{equation}\label{eqn.xkyk}
\mathrm{rank}\left(X^k\right) -\mathrm{rank}\left(J_\mu^k\right)
=\begin{cases}
0, &\hbox{if $v_{\alpha_i}=0$ for all $i$ such that $\mu_i \geq k$},\\
1, &\hbox{otherwise.}
\end{cases}
\end{equation}
Let $c \leq n$ be the smallest positive integer 
for which $\mathrm{rank}(X^c)-\mathrm{rank}(J_\mu^c)=0$.  Then Equation~\eqref{eqn.xkyk} implies that
$$\mathrm{rank}(X^k)-\mathrm{rank}(J_\mu^k)=
\begin{cases}
0, &\hbox{if $k \geq c$},\\
1, &\hbox{if $k<c$}.
\end{cases}
$$
Together with Lemma~\ref{lem.rank}, we deduce that
$$
\lambda_k'-\mu_k'
= \left(\mathrm{rank}(X^{k-1}) - \mathrm{rank}(X^{k}) \right)
	- \left(\mathrm{rank}(J_\mu^{k-1}) - \mathrm{rank}(J_\mu^{k}) \right)
= \begin{cases}
1, &\hbox{if $k=c$},\\
0, &\hbox{if $k\neq c$}.
\end{cases}
$$
Therefore, $\lambda \gtrdot \mu$ in the case $Y=J_\mu$.

In the general case where $Y$ is any matrix of Jordan type $\mu$, then $\mathrm{rank}(Y^k) = \mathrm{rank}(J_\mu^k)$ for all $k\geq0$, so the argument is the same.
\qed


Let $\lambda$ be the partition whose diagram is obtained by adding a box to the $i$th row and $j$th column of the diagram of the partition $\mu$.
Define the coefficient 
\begin{equation}\label{eq.c}
c_{\mu,\lambda}(q) 
= \left\{ 
\begin{array}{ll}
q^{|\mu|-\mu_j'}, &\hbox{if $j=1$},\\
q^{|\mu|-\mu_{j-1}'}\left(q^{\mu_{j-1}'-\mu_j'}-1 \right),
&\hbox{if $j\geq2$}.
\end{array}
\right.
\end{equation}
Note that in the case $j\geq2$, we have $\mu_{j-1}'-\mu_j' \geq1$.

\begin{lemma}\label{lem.2} 
Let $Y$ be an upper-triangular nilpotent matrix of Jordan type $\mu\vdash n-1$.  If $\mu\lessdot\lambda$, then there are $c_{\mu,\lambda}(q)$ upper-triangular nilpotent matrices $X$ of Jordan type $\lambda$ whose first order leading principal submatrix is $Y$.
\end{lemma}
\proof 

By similarity, it suffices to consider the case $Y=J_\mu = J_{\mu_1}\oplus \cdots \oplus J_{\mu_m}$, where $\ell(\mu)=m$.
Suppose $X$ is a matrix of the form
$$X= 
\begin{bmatrix} 
J_\mu & \mathbf{v}\\
\mathbf{0} &0
\end{bmatrix}
$$
of Jordan type $\lambda$ such that $\lambda$ is obtained by adding a box to $\mu$ in the $i$th row and $j$th column.


First consider the case $j\geq2$.  
Following the proof of Lemma~\ref{lem.1}, we know that $j$ is the unique integer where
$\mathrm{rank}(X^{j-1}) = \mathrm{rank}(J_\mu^{j-1})+1,$ and
$\mathrm{rank}(X^{j}) = \mathrm{rank}(J_\mu^{j})$.
In order to satisfy the first condition, the entries in the vector $\mathbf{v}^{j-1}$ corresponding to the boxes in the $(j-1)$th column and rows $\geq i$ must not simultaneously be zero (refer to Equation~\eqref{eqn.visualizev} and Example~\ref{eg.visualize}), while in order to satisfy the second condition, the entries in the vector $\mathbf{v}^j$ corresponding to the boxes in the $j$th column of $\mu$ must all be zero.  The remaining $n-1-\mu_{j-1}'$ entries of the vector $\mathbf{v}$ are free to be any element in $\mathbb{F}_q$, so there are
$$q^{n-1-\mu_{j-1}'} \left(q^{\mu_{j-1}'-\mu_j'} -1\right) $$
possible matrices $X$ whose leading principal submatrix is $J_\mu$.



The case $j=1$ is simpler.  
The necessary and sufficient condition that $X$ and $J_\mu$ must satisfy is that $\mathrm{rank}(X^k) = \mathrm{rank}(J_\mu^k)$ for all $k\geq1$, so the entries in the vector $\mathbf{v}$ corresponding to the boxes in the first column of the diagram for $\mathbf{v}^1$ must all be zero, while the remaining $n-1-\mu_1'$ entries are free to be any element in $\mathbb{F}_q$, so there are $q^{n-1-\mu_1'}$ matrices $X$ whose leading principal submatrix is $J_\mu$ in this case.
\qed

\begin{theorem}\label{thm.JFrec}
The number of $n$ by $n$ upper-triangular nilpotent matrices over $\mathbb{F}_q$ of Jordan type $\lambda\vdash n$ is
$$F_\lambda(q) = \sum_{\mu:\, \mu\lessdot
\lambda} c_{\mu,\lambda}(q) F_\mu(q), $$
with $F_{\emptyset}(q) =1$.
\end{theorem}
\proof
Proceed by induction on $n$.  For $n=1$, the zero matrix is the only upper-triangular nilpotent matrix, and it has Jordan type $(1)$, agreeing with the formula $c_{\emptyset,(1)}(q) = 1$.  

Suppose $\lambda\vdash n$. By Lemma~\ref{lem.1}, any matrix of Jordan type $\lambda$ has a leading principal submatrix of type $\mu\vdash n-1$ for some $\mu\lessdot\lambda$.  Furthermore, by Lemma~\ref{lem.2}, for each matrix $Y\in \mathfrak{g}_{n-1,\mu}(\mathbb{F}_q)$, there are $c_{\mu\lambda}(q)$ matrices $X\in\mathfrak{g}_{n,\lambda}(\mathbb{F}_q)$ having $Y$ as its leading principal submatrix.  Summing over all $\mu \lessdot\lambda$ gives the desired formula.
\qed

\begin{figure}
\Yboxdim8pt
$$\xymatrixcolsep{4pc}\xymatrixrowsep{3pc}
\xymatrix{
& {\yng(1,1,1,1)}\\
& {\yng(1,1,1)} \ar[r]^{q^3-1}		\ar[u]^-1
	& {\yng(2,1,1)} 
	& {\yng(2,2)} \\
& {\yng(1,1)} \ar[r]^{q^2-1}		\ar[u]^-1
	& {\yng(2,1)} \ar[r]^{(q-1)q^2} \ar[u]^-q 		\ar[ur]|{(q-1)q}
	& {\yng(3,1)} \\
& {\yng(1)} \ar[r]^{q-1} 			\ar[u]^-1
	& {\yng(2)} \ar[r]^{(q-1)q} 	\ar[u]^-q
	& {\yng(3)} \ar[r]^{(q-1)q^2} 	\ar[u]^-{q^2}
	& {\yng(4)} \\
\emptyset \ar[ur]^1 
}
$$
\caption{Young's lattice with edge weights $c_{\mu,\lambda}(q)$, up to $n=4$.}
\label{fig.Y}
\end{figure}
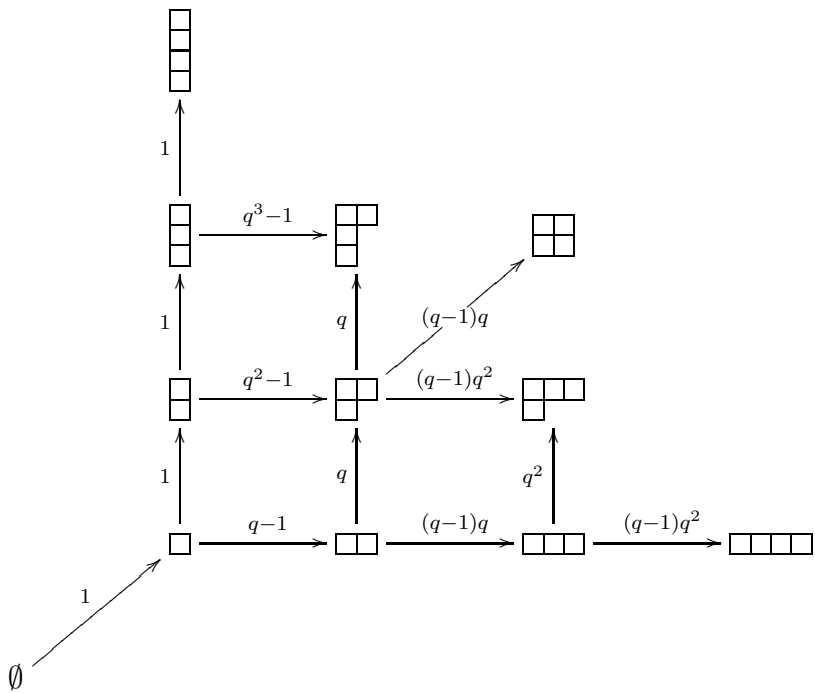

\begin{remark}[Formulation in terms of standard Young tableaux]

The formula for $F_\lambda(q)$ in Theorem~\ref{thm.JFrec} can be re-phrased as a sum over the set $\mathcal{P}_{\mathcal{Y}}(\lambda)$ of paths in Young's lattice $\mathcal{Y}$ from the empty partition $\emptyset$ to $\lambda$.  
If $\mu \lessdot\lambda$ in $\mathcal{Y}$, we assign the weight 
$c_{\mu,\lambda}(q)$ to the corresponding edge in $\mathcal{Y}$.  
Figure~\ref{fig.Y} shows Young's lattice with weighted edges for partitions with up to four boxes.
Let $P = (\emptyset = \pi^{(0)}, \pi^{(1)}, \ldots, \pi^{(n)}=\lambda)$ denote a path in $\mathcal{Y}$ from $\emptyset$ to $\lambda$, where $\pi^{(i)}$ is a partition of $i$.  To simplify notation, let $\epsilon_i(q)= c_{\pi^{(i-1)}, \pi^{(i)}}(q)$. 
Theorem~\ref{thm.JFrec} is equivalently re-phrased as
\begin{equation}\label{eq.paths}
F_\lambda(q) = \sum_{P\in P_{\mathcal{Y}}(\lambda)} F_P(q),
\end{equation}
where the weight of the path $P$ is $F_P(q) = \prod_{i=1}^n \epsilon_i(q)$.

The set of paths $\mathcal{P}_{\mathcal{Y}}(\lambda)$ is in bijection with the set $\mathrm{SYT}(\lambda)$ of standard Young tableaux of shape $\lambda$, so we can also give an equation for $F_\lambda(q)$ as a sum over standard Young tableaux.

A {\em standard Young tableau} $T$ of shape $\lambda$ is a filling of the Ferrers diagram of $\lambda\vdash n$ with the integers $1,\ldots, n$ such that the integers increase weakly along each row and strictly along each column. For $1\leq i\leq n$,  Let $T^{(i)}$ denote the Young tableau of shape $\lambda^{(i)}$ consisting of the boxes containing $1,\ldots, i$, and define weights
\begin{equation}
T^{(i)}(q) =
\begin{cases}
q^{i-\ell(\lambda^{(i)})}, 
	& \hbox{if the $i$th box is in the first column,}\\
q^{i-{\lambda_{j}^{(i)}}'} - q^{i-1-{\lambda_{j-1}^{(i)}}'},
 	& \hbox{if the $i$th box is in the $j$th column, $j\geq2$.}
\end{cases}
\end{equation}
Then
\begin{equation}\label{eqn.syt}
F_\lambda(q)= \sum_{T\in \mathrm{SYT}(\lambda)} F_T(q),
\end{equation}
where the weight of the standard Young tableau $T$ is $F_T(q) = \prod_{i=1}^n T^{(i)}(q)$.
\end{remark}


\subsection{Properties of $F_\lambda(q)$}
Several properties of $F_\lambda(q)$ follow readily from Theorem~\ref{thm.JFrec}.  
For $\lambda \vdash n$, let
\begin{equation}\label{eq.nlambda}
n_\lambda= \sum_{i\geq1} (i-1)\lambda_i = \sum_{b\in \lambda} \mathrm{coleg}(b),
\end{equation}
where if a box $b\in \lambda$ lies in the $i$th row of $\lambda$, then $\mathrm{coleg}(b) = i-1$.

\begin{corollary}\label{cor.degree}
Let $\lambda \vdash n$.  
As a polynomial in $q$, 
$$\deg F_\lambda(q) =  {n\choose 2}- n_\lambda.$$
Moreover, the coefficient of the highest degree term in $F_\lambda(q)$ is $f^\lambda$, the number of standard Young tableaux of shape $\lambda$.
\end{corollary}
\proof
Suppose $P = (\emptyset = \pi^{(0)}, \pi^{(1)}, \ldots, \pi^{(n)}=\lambda)$ is a path in $\mathcal{Y}$ such that $\pi^{(k)}$ is obtained by adding a box to the $i$th row and $j$th column of $\pi^{(k-1)}$.  Then $\deg c_{\pi^{(k-1)},\pi^{(k)}}(q)= k-i$, and therefore
$$\deg F_P(q) = \sum_{k=1}^n \deg c_{\pi^{(k-1)},\pi^{(k)}}(q)
= \sum_{k=1}^n k - \sum_{k\geq1}k\lambda_k
= {n\choose 2} - n_\lambda.$$
In particular, every polynomial $F_P(q)$ arising from a path $P \in \mathcal{P}_{\mathcal{Y}}(\lambda)$ has the same degree, so $\deg F_\lambda(q) =  {n \choose 2} - n_\lambda$.  Moreover, each $F_P(q)$ is monic, so the coefficient of the highest degree term in $F_\lambda(q)$ is the number of paths in $\mathcal{Y}$ from $\emptyset$ to $\lambda$, which is $f^\lambda$.
\qed

\begin{corollary}\label{cor.multiplicity}
Let $\lambda \vdash n$.  
The multiplicity of the factor $q-1$ in $F_\lambda(q)$ is $n-\ell(\lambda)$.
\end{corollary}
\proof
The weight $c_{\pi^{(k-1)},\pi^{(k)}}(q)$ corresponding to the $k$th step in the path $P$ contributes a single factor of $q-1$ to $F_P(q)$ if and only if the $k$th box added is not in the first column of $\lambda$.
Therefore, the multiplicity of $q-1$ in $F_P(q)$ is $n-\ell(\lambda)$, and it follows that the multiplicity of $q-1$ in $F_\lambda(q)$ is $n-\ell(\lambda)$.
\qed

\begin{example}
There are two partitions of $4$ with two parts, namely $(3,1)$ and $(2,2)$.

There are three paths from $\emptyset$ to $(3,1)$ in $\mathcal{Y}$, giving
\begin{align*}
F_{(3,1)}(q) 
&= (q-1) \cdot (q-1)q \cdot q^2
+ (q-1) \cdot q \cdot (q-1)q^2
+ \cdot (q^2-1) \cdot (q-1)q^2 \\
&= (q-1)^2\left(3q^3+q^2\right),
\end{align*}
and there are two paths from $\emptyset$ to $(2,2)$ in $\mathcal{Y}$, giving
\begin{align*}
F_{(2,2)}(q) &= (q-1)\cdot q\cdot (q-1)q
+ (q^2-1)\cdot (q-1)q \\
&= (q-1)^2(2q^2+q).
\end{align*}
Summing these gives a shift of the $q$-Stirling polynomial
$(q-1)^2q^4S_{4,2}(q^{-1})=(q-1)^2(3q^3+3q^2+q)$.
\end{example}

\subsection{Explicit formulas}
In this section, we derive non-recursive formulas for some special cases of $\lambda$.  Previously, we have noted the simple cases $F_{(1^n)}= 1$ and $F_{(n)} = q^{n-1\choose 2}(q-1)^{n-1}$.

\begin{prop}[Hook shapes] \label{prop.hookshapes}
Let $n> k \geq2$, and let $\lambda = (n-k+1, 1^{k-1})$ be a hook-shaped partition of $n$ with $\ell(\lambda) = k$ parts.  
Then
$$F_\lambda(q) 
= (q-1)^{n-k} \sum_{i=0}^{k-1} {n-i-1\choose k-i-1} q^{\alpha -i},
\qquad
\hbox{where }
\alpha = {n-1\choose 2}-{{k-1} \choose 2}.
$$
\end{prop}
\proof
We make use of Equation~\eqref{eq.paths}.  
We enumerate paths from $\emptyset$ to $\lambda$ according to the first time a box is added to the second column,
so for $0\leq r\leq k-1$, let $\mathrm{S}_r$ be the set of paths in the sublattice $[\emptyset,\lambda]$ which contains the edge $((1,1^r), (2,1^r))$.  Such paths are formed by the concatenation of the unique path between $\emptyset$ and $\nu= (2,1^r)$, which has weight $q^{r+1}-1$, with any path in the sublattice $[\nu,\lambda]$.
The sublattice $[\nu,\lambda]$ is the Cartesian product of a $(n-k)$-chain and a $(k-r-1)$-chain, so it forms a rectangular grid, and therefore 
$|\mathrm{S}_r|= {n-r-1 \choose k-r-1}$.
Notice that in any sublattice of the form
$$\xymatrixcolsep{4pc}\xymatrixrowsep{3pc}
\xymatrix{
(a,1^{b+1}) 					\ar[r]^{(q-1)q^{a+b}} 
	& (a+1,1^{b+1})
	\\
(a,1^b) \ar[r]^{(q-1)q^{a+b-1}} \ar[u]_-{q^{a-1}}
	& (a+1,1^b)					\ar[u]^-{q^{a}}
}
$$
the product of the edge weights is $(q-1)q^{2a+b-1}$ no matter which path is taken from $(a,1^b)$ to $(a+1,1^{b+1})$, so it follows that every path from $\nu$ to $\lambda$ has the same weight. 
By considering the path $(\nu, 21^{r+1}, \ldots, 21^{k-1}, 31^{k-1},\ldots, \lambda)$, this weight is easily seen to be
$(q-1)^{n-k-1} q^{\alpha-r},$ for 
$\alpha={n-1\choose 2} -{k-1\choose 2}$.
Altogether, 
$$F_\lambda(q)
= (q-1)^{n-k} \sum_{r=0}^{k-1} {n-r-1\choose k-r-1} 
\left(q^{\alpha}+q^{\alpha-1}+\cdots+q^{\alpha-r}\right).$$
For $0 \leq i \leq k-1$, the coefficient of $q^{\alpha-i}$ in $F_\lambda(q)/(q-1)^{n-k}$ is 
$$\sum_{r=i}^{k-1} {n-r-1 \choose k-r-1} 
= \sum_{u=0}^{k-i-1} {n-k+u\choose u} 
= {n-i-1\choose k-i-1},$$
since $\sum_{u=0}^M {N+u\choose u} = {N+M+1\choose M}$.
Therefore,
$$F_\lambda(q) =
(q-1)^{n-k} \sum_{i=0}^{k-1} {n-i-1\choose k-i-1} q^{\alpha-i},$$
as claimed.
\qed
 
We next consider the case when $\lambda$ is a partition with two parts.  For $n \geq k \geq 1$, let
\begin{equation}\label{eqn.gencatalan}
C_{n,k} = {n+k\choose k} - {n+k \choose k-1},
\end{equation}
and let $C_{n,0} =1$ for all $n\geq0$. These generalized Catalan numbers $C_{n,k}$ enumerate lattice paths from $(0,0)$ to $(n,k)$, using the steps $(1,0)$ and $(0,1)$, which do not rise above the line $y=x$.  
In the remainder of this section, we shall refer to these as {\em Dyck paths}.  

The generalized Catalan numbers satisfy the simple recursive formula $C_{n,k} = C_{n,k-1}+C_{n-1,k}.$  
Also, these are the usual Catalan numbers
$C_n = \frac{1}{n+1}{2n\choose n}=C_{n,n}= C_{n,n-1}$
when $k=n$ or $n-1$.
These facts will be used in the computations which follow.

\begin{figure}
$$
%
\begin{array}{c|ccccccc}
n\backslash k &0&1&2&3&4&5&6\\
\hline
0&1&&&&&&\\
1&1&1&&&&&\\
2&1&2&2&&&&\\
3&1&3&5&5&&&\\
4&1&4&9&14&14&&\\
5&1&5&14&28&42&42&\\
6&1&6&20&48&90&132&132
\end{array}
$$
\caption{The Catalan triangle $C_{n,k}$.}
\label{fig.catalantriangle}
\end{figure}

\begin{prop}[Partitions with two parts]
If $\lambda = (r,s) \vdash n$ such that $r>s\geq1$, then
$$F_{(r,s)}(q) = (q-1)^{r+s-2} q^{{r+s-1\choose 2}-2s+1}
\sum_{i=0}^s C_{r+s-i,i}q^i.
$$
If $r=s$, then
$$F_{(r,r)}(q) = (q-1)^{2r-2} q^{{2r-2\choose 2}} 
	\sum_{i=0}^{r-1} C_{2r-1-i,i}q^i.
$$
\end{prop}
\proof Proceed by induction on $r+s$.  
The base cases are $F_{(r)}(q) = (q-1)^{r-1}q^{r-1\choose 2}$ for $r\geq1$ and $F_{(1,1)}(q) = 1$.

We first handle the case $s=1$ separately. For $r\geq2$,
\begin{align*}
F_{(r,1)}(q)
&= q^{r-1}F_{(r)}(q)+ (q-1)q^{r-1}F_{(r-1,1)}(q)\\
&= q^{r-1} \cdot (q-1)^{r-1}q^{r-1\choose 2} 
	+ (q-1)q^{r-1}\cdot (q-1)^{r-2} q^{{r-1\choose2}-1} 
		\left( (r-1)q+ 1 \right)\\
&= (q-1)^{r-1}q^{{r\choose2}-1} \left( rq +1 \right).
\end{align*}

Next, consider the case $s=r$.  For $r\geq2$,
\begin{align*}
F_{(r,r)}(q)
&= (q-1) q^{2r-3}F_{(r,r-1)}(q)\\
&= (q-1) q^{2r-3}\cdot (q-1)^{2r-3} q^{{2r-2\choose 2}- 2(r-1) +1} 
	\sum_{i=0}^{r-1} C_{2r-1-i,i}q^i\\
&= (q-1)^{2r-2} q^{{2r-2\choose 2}} 
	\sum_{i=0}^{r-1} C_{2r-1-i,i}q^i.
\end{align*}

The case $s=r-1$ is obtained as follows.  For $r\geq3$,
\begin{align*}
F_{(r,r-1)}(q)
&= (q-1)q^{2r-4}F_{(r,r-2)}(q) + (q^2-1) q^{2r-4} F_{(r-1,r-1)}(q)\\
&= (q-1)^{2r-3}q^{2r-4}q^{2r-4\choose2} 
	\left(q\sum_{i=0}^{r-2} C_{2r-2-i,i} q^i 
		+ (q+1)\sum_{i=0}^{r-2} C_{2r-3-i,i}q^i \right)\\
&= (q-1)^{2r-3}q^{2r-3\choose2} 
	\Bigg( \left(C_{r,r-2}+C_{r-1,r-2} \right) q^{r-1}\\
 &\qquad 
 + \sum_{i=1}^{r-2}
 	\left(C_{2r-1-i,i-1}+C_{2r-2-i,i-1}+C_{2r-3-i,i} \right) q^i
	+ C_{2r-3,0}q^0\Bigg).
\end{align*}
Since $C_{n,n-1}=C_{n,n}$, then
$C_{r,r-2}+C_{r-1,r-2} = C_{r,r-1}$.  Similarly, we obtain
$C_{2r-1-i,i-1}+C_{2r-2-i,i-1}+C_{2r-3-i,i} =  C_{2r-1-i,i}$ by
applying the recurrence equation for the generalized Catalan numbers.  Lastly, $C_{n,0}=1$ for all $n\geq0$, thus
\begin{align*}
F_{(r,r-1)}(q)
&= (q-1)^{2r-3}q^{2r-3\choose2} 
	\sum_{i=0}^{r-1}C_{2r-1-i,i} q^i,
\end{align*}
which agrees with the formula for the case $s=r-1$.

The last case to consider is the general case $r-s\geq2$ where $s\geq2$.
\begin{align*}
F_{(r,s)}
&= (q-1)q^{r+s-3} F_{(r,s-1)} + (q-1) q^{r+s-2} F_{(r-1,s)}\\
&=  (q-1)^{r+s-2}q^{r+s-3} q^{{r+s-2\choose2}-2s+1}
\left(q^2 \sum_{i=0}^{s-1} C_{r+s-1-i,i}q^i
	+ q\sum_{i=0}^s C_{r-1+s-i,i}q^i\right)\\
&= (q-1)^{r+s-2}q^{{r+s-1\choose2}-2s+1}
	\left(C_{r+s,0}q^0+\sum_{i=1}^{s} 
		\left(C_{r+s-i,i-1}+C_{r-1+s-i,i}\right)q^i
	\right)\\
&= (q-1)^{r+s-2}q^{{r+s-1\choose2}-2s+1}
	\sum_{i=0}^{s} C_{r+s-i,i} q^i.
\end{align*}
\qed


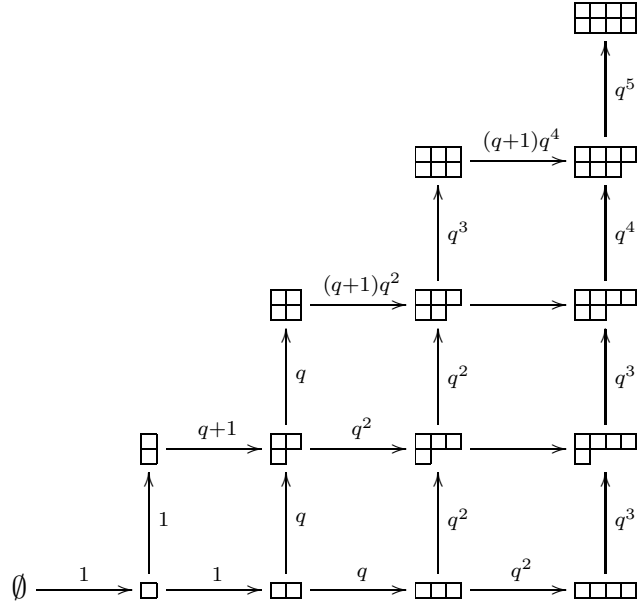
\begin{figure}
\Yboxdim6pt
$$
\xymatrixcolsep{3pc}\xymatrixrowsep{3pc}
\xymatrix{
&&&& {\yng(4,4)}
	\\
&&&	{\yng(3,3)} 	\ar[r]^{(q+1)q^4}
	& {\yng(4,3)}					\ar[u]_{q^5}
	\\
&&	{\yng(2,2)} 	\ar[r]^{(q+1)q^2}
	& {\yng(3,2)}	\ar[r]^{}		\ar[u]_{q^3}
	& {\yng(4,2)}					\ar[u]_{q^4}
	\\
&	{\yng(1,1)}		\ar[r]^{q+1}
	& {\yng(2,1)}	\ar[r]^{q^2}	\ar[u]_{q}
	& {\yng(3,1)}	\ar[r]			\ar[u]_{q^2}
	& {\yng(4,1)}					\ar[u]_{q^3}
	\\
\emptyset 			\ar[r]^1
	& {\yng(1)}		\ar[r]^{1}		\ar[u]_1
	& {\yng(2)}		\ar[r]^{q}		\ar[u]_q
	& {\yng(3)}		\ar[r]^{q^2}	\ar[u]_{q^2}
	& {\yng(4)}						\ar[u]_{q^3}
}
$$\caption{Factors of $q-1$ are omitted from the edge weights in this sublattice of partitions with at most two rows.}
\label{fig.tworows}
\end{figure}

The next equation is a formula for $F_{(r,r)}(q)$ with a different factorization.
\begin{prop}[Two equal parts] \label{prop.kewkewplusone}
Let $\lambda = (r,r)\vdash n$, and $r\geq1$.
Then
$$F_{(r,r)}(q)
= (q-1)^{2r-2} 
	\sum_{i=0}^{r-1} C_{r-1,r-1-i}\ q^{2(r-1)^2-i}(q+1)^i.
$$
\end{prop}
\proof The set of paths in the sublattice $[\emptyset,\lambda]$ are in bijection with the set of lattice paths from $(0,0)$ to $(r,r)$.
In any sublattice of the form
$$\xymatrixcolsep{4pc}\xymatrixrowsep{3pc}
\xymatrix{
(a,b+1) 							\ar[r]^{(q-1)q^{a+b}} 
	& (a+1,b+1)
	\\
(a,b) \ar[r]^{(q-1)q^{a+b-1}} 		\ar[u]^-{(q-1)q^{a+b-2}}
	& (a+1,b)						\ar[u]_-{(q-1)q^{a+b-1}}
}
$$
where $b\leq a-2$, the product of the edge weights is $(q-1)^2q^{2a+2b-2}$ no matter which path is taken from $(a,b)$ to $(a+1,b+1)$.  As for sublattices of the form
$$\xymatrixcolsep{4pc}\xymatrixrowsep{3pc}
\xymatrix{
(a,a) 							\ar[r]^{(q^2-1)q^{2a-2}} 
	& (a+1,a)\\
(a,a-1) 			\ar[r]^{(q-1)q^{2a-2}} 	\ar[u]^-{(q-1)q^{2a-3}}
	& (a+1,a-1)					\ar[u]_-{(q-1)q^{2a-2}}
}
$$
the product of the edge weights is $(q-1)^2 q^{4a-4}$ via the lower horizontal edge, versus $(q-1)^2 q^{4a-5} (q+1)$ via the upper horizontal edge.  It follows that if a path $P$ from $\emptyset$ to $\lambda$ contains $i$ partitions of the form $(a,a)$, then it has the weight
$$F_p(q)  = (q-1)^{2r-2} q^{2(r-1)^2-i} (q+1)^i.$$

Dyck paths may be enumerated according to the points at which they touch the diagonal line $y=x$, and the set of touch points are indexed by compositions $\alpha= (\alpha_1,\ldots, \alpha_{i+1})\vDash r$ where $\alpha_j \geq 1$.
The number of Dyck paths from $(0,0)$ to $(r,r)$ which touch the diagonal exactly $i$ times, not including the initial and the end points, is
$$\sum_{\alpha \vDash r \atop \ell(\alpha)=i+1} \prod_{j=1}^{i+1} C_{\alpha_j-1}.$$
On the other hand, the number $C_{r-1, r-1-i}$ of Dyck paths from $(0,0)$ to $(r-1, r-1-i)$ satisfies the same recurrence equation
$$C_{r-1, r-1-i} = \sum_{\beta \vDash r-1-i\atop \ell(\beta) = i+1}
\prod_{j=1}^{i+1} C_{\beta_j},$$
but the sum is over the set of weak compositions so that $\beta_j \geq0$.  Under the appropriate shift in indices, it follows that the number of Dyck paths from $(0,0)$ to $(r,r)$ which touch the diagonal exactly $i$ times is $C_{r-1,r-1-i}$.
The result follows from this.\qed

\begin{corollary}
For $k\geq m \geq0$,
$$\sum_{j=m}^k {j\choose m} C_{k,j}=C_{2k+1-m,m}. $$
\end{corollary}
\proof The two formulas for $F_{(k+1,k+1)}(q)$ yields the identity
$$\sum_{i=0}^{k} C_{k,k-i}\ q^{2k^2-i}(q+1)^i
=q^{{2k\choose 2}} \sum_{i=0}^{k} C_{2k-i,i}q^i.
$$
Extracting the coefficient of $q^{2k^2-m}$ in the above expressions yields the result.\qed

\begin{remark} The formula for $F_{(r,r)}(q)$ provided in
Proposition~\ref{prop.kewkewplusone} can be viewed as a sum over Dyck paths, where each Dyck path $\pi$ contributes a term of the form $q^{s_1(\pi)} (q+1)^{s_2(\pi)}$ for some statistics $s_1$ and $s_2$ on the Dyck paths.
This particular factorization for $F_{(r,r)}(q)$ is related to the work of Cai and Readdy on the $q$-Stirling numbers of the second kind, since the polynomials $F_\lambda(q)$ can be viewed as a refinement of $S_{n,k}(q)$, as explained in Section~\ref{sec.rooks}.

Cai and Readdy obtained a formula~\cite[Theorem 3.2]{caiReaddy15} for $\widetilde S_{n,k}(q)$ (they use a different recursive formula to define the $q$-Stirling numbers, and the two are related by $S_{n,k}(q) = q^{k\choose 2}\widetilde S_{n,k}(q)$) as a sum over allowable restricted-growth words, where each allowable word $w$ gives rise to a term of the form $q^{a(w)}(q+1)^{b(w)}$ for some statistics $a(w)$ and $b(w)$.  They also showed that this enumerative result has an interesting extension to the study of the Stirling poset of the second kind, providing a decomposition of that poset into Boolean sublattices.

For example, if we let $(q-1)^{n-\ell(\lambda)}G_\lambda(q) = F_\lambda(q)$, then
$G_{(3,1)}(q) + G_{(2,2)}(q)  = q^3 \widetilde S_{4,2}(q^{-1}).$
The formula of Cai and Readdy yields $q^3 S_{4,2}(q^{-1}) = q(q+1)^2 + q^2(q+1)+q^3$, while our factorization yields 
$G_{(3,1)}(q) + G_{(2,2)}(q) 
= \left(q^3+q^3+q^2(q+1) \right) + \left(q^2+ q(q+1)\right)$.
So, the result of Proposition~\ref{prop.kewkewplusone} gives a different expression for $S_{n,n-2}(q)$ as a sum 
with terms of the form $q^{s_1(\pi)} (q+1)^{s_2(\pi)}$, and it may be interesting to further investigate such factorizations of $F_\lambda(q)$.
\end{remark}

\begin{example}  The first few $F_{(k,k)}$ are
\begin{align*}
F_{(1,1)} &= 1\\
F_{(2,2)} 
	&= (q-1)^2\left(q^2+ q(q+1)\right)\\
	&= (q-1)^2\left(2q^2+q\right)\\ 
F_{(3,3)} 
	&= (q-1)^4\left(2q^8+ 2q^7(q+1)+ q^6(q+1)^2\right)\\
	&= (q-1)^4\left(5q^8+ 4q^7+ q^6\right)\\
F_{(4,4)} 
	&= (q-1)^6\left(5q^{18}+ 5q^{17}(q+1)+ 3q^{16}(q+1)^2+ q^{15}(q+1)^3\right)\\
	&= (q-1)^6\left(14q^{18}+ 14q^{17}+ 6q^{16}+ q^{15}\right)\\
F_{(5,5)} 
	&= (q-1)^8\left(14q^{32}+ 14q^{31}(q+1)+ 9q^{30}(q+1)^2+ 4q^{29}(q+1)^3+ q^{28}(q+1)^4\right)\\
	&= (q-1)^8\left(42q^{32}+ 48q^{31}+ 27q^{30}+ 8q^{29}+ q^{28}\right).
\end{align*}
\end{example}

We end this section with one more closed formula for $F_\lambda(q)$ where $\lambda$ is a rectangular shape with two columns. Let $\mathcal{D}(n,k)$ denote the set of Dyck paths from $(0,0)$ to $(n,k)$.  The {\em coarea} of a Dyck path $\pi$ is the number of whole unit squares lying between the path and the $x$-axis.  For $i=1,\ldots, n$, let $\rho_i(\pi)$ be one plus the number of unit squares lying between the path and the line $y=x+1$ in the $i$th row.
For example, the following Dyck path $\pi$ has $\mathrm{coarea}(\pi)=12$, and $(\rho_1(\pi),\rho_2(\pi),\rho_3(\pi),\rho_4(\pi)) = (2,2,1,1)$.
$$\beginpicture
\setcoordinatesystem units <.7cm,.7cm> 
\setplotarea x from -1.5 to 6, y from -.5 to 4.5
\plot 0 0 6 0 / \plot 0 1 6 1 / \plot 0 2 6 2 / \plot 0 3 6 3 /
\plot 0 4 6 4 / 
\plot 0 0 0 4 / \plot 1 0 1 4 / \plot 2 0 2 4 / \plot 3 0 3 4 /
\plot 4 0 4 4 / \plot 5 0 5 4 / \plot 6 0 6 4 /
\put{\tiny$(0,0)$} at 0 -.5
\put{\tiny$y=x+1$} at -1.5 .5
{\setdashes
\plot 0 0 4 4 /
\plot -.5 .5 3.5 4.5 /
}\setplotsymbol(.)
\plot 0 0 2 0 2 1 3 1 3 3 4 3 4 4 6 4 /
\endpicture$$
For $n\geq1$, let $[n]_q = 1+q+\cdots +q^{n-1}$.

\begin{figure}
\Yboxdim6pt
$$
\xymatrixcolsep{2.5pc}\xymatrixrowsep{1.5pc}
\xymatrix{
&&&& {\yng(2,2,2,2)}
	\\
&&& {\yng(2,2,2)} 		\ar[r]^{q^3}
	& {\yng(2,2,2,1)}					\ar[u]_{q^3[1]}
	\\
&& {\yng(2,2)}			\ar[r]^{q^2}
	& {\yng(2,2,1)} 	\ar[r]^{q^2}	\ar[u]_{q^2[1]}
	& {\yng(2,2,1,1)}					\ar[u]_{q^2[2]}
	\\
& {\yng(2)}			\ar[r]^q
	& {\yng(2,1)}	\ar[r]^q 	\ar[u]_{q[1]}
	& {\yng(2,1,1)} \ar[r]^q	\ar[u]_{q[2]}
	& {\yng(2,1,1,1)}			\ar[u]_{q[3]}
	\\
\emptyset 			\ar[r]^1
	& {\yng(1)}		\ar[r]^1	\ar[u]_{[1]}
	& {\yng(1,1)}	\ar[r]^1	\ar[u]_{[2]}
	& {\yng(1,1,1)}	\ar[r]^1	\ar[u]_{[3]}
	& {\yng(1,1,1,1)}			\ar[u]_{[4]}
}
$$
\caption{Factors of $q-1$ are omitted from the edge weights in this sublattice of partitions with at most two columns.}
\label{fig.twocolumns}
\end{figure}

\begin{prop}(Partitions with two columns) Let $\lambda = (2^r, 1^s) \vdash n$ such that $r,s \geq 0$.  Then
$$F_{(2^r,1^s)}(q) = (q-1)^r q^{r\choose2} 
\sum_{\pi\in \mathcal{D}(r+s,r)}
q^{\mathrm{coarea}(\pi)} \prod_{i=1}^r \left[\rho_i(\pi)\right]_q.$$
\end{prop}

\proof By Corollary~\ref{cor.multiplicity}, we know the multiplicity of the factor $q-1$ in $F_{(2^r,1^s)}(q)$ is $n-\ell(\lambda) = \lambda_2' = r$, so we focus on computing $F_{(2^r,1^s)}(q)/(q-1)^r$.  The paths in Young's lattice from $\emptyset$ to $(2^r,1^s)$ are in bijection with the Dyck paths $\mathcal{D}(r+s,r)$, so we identify these paths; adding a box in the first column of a partition corresponds to a $(1,0)$ step in the Dyck path, and adding a box in the second column of a partition corresponds to a $(0,1)$ step in the Dyck path.  As seen in Figure~\ref{fig.twocolumns}, a vertical step $(i,j)$ to $(i,j+1)$ has weight $q^j[i]_q$, while a horizontal step $(i,j)$ to $(i+1,j)$ has weight $q^j$.  Thus the product of the edge weights of the $r$ vertical steps of a given Dyck path $\pi$ is $q^{r\choose2}\prod_{i=1}^r [\rho_i(\pi)]_q$, while the product of the edge weights of the $r+s$ horizontal steps of a given Dyck path is $q^{\mathrm{coarea}(\pi)}$.  The result follows.
\qed

\begin{example}
The first few $F_{(2^n)}$ are
\begin{align*}
F_{(2)} &= (q-1)(q+1)\\
F_{(2^2)} &= (q-1)^2 q \left(q+ (q+1) \right)\\
F_{(2^3)} &= (q-1)^3 q^3 
	\left(q^3+ 2q^2(q+1)+ q(q+1)^2+ (q^2+q+1)(q+1)\right)\\
F_{(2^4)} &= (q-1)^4 q^6 
	\left(q^6+ 3q^5[2] + 3q^4[2]^2 + q^3[2]^3 + 2q^3[3]! 
	+ 2q^2[2][3]! + q[3][3]!+ [4]!\right).
\end{align*}
\end{example}

\begin{remark}
Kirillov and Melnkov~\cite{kirillovMelnikov95} considered 
the number $A_n(q)$ of $n$ by $n$ upper-triangular matrices over $\mathbb{F}_q$ satisfying $X^2=0$.  In their first characterization of these polynomials, they considered the number $A_n^r(q)$ of matrices of a given rank $r$, so that $A_n(q) = \sum_{r\geq0}A_n^r(q)$, and observed that $A_n^r(q)$ satisfies the recurrence equation
\begin{displaymath}
A_{n}^{r}(q) = q^r A_{n-1}^r(q) +\left(q^{n-r}-q^r\right) A_n^r(q), \qquad
A_n^0(q)=1.
\end{displaymath}
We may think of $A_n(q)$ as the sum of $F_\lambda(q)$ over $\lambda \vdash n$ with at most two columns, so Theorem~\ref{thm.JFrec} is a generalization of this recurrence equation.

It was also conjectured in~\cite{kirillovMelnikov95} that the same sequence of polynomials arise in a number of different ways. Ekhad and Zeilberger~\cite{ekhadZeilberger96} proved that one of the conjectured alternate definitions of $A_n(q)$, namely
\begin{displaymath}
C_n(q) = \sum_{s} c_{n+1,s}q^{\frac{n^2}{4}+\frac{1-s^2}{12}},
\end{displaymath}
is a sum over all $s\in [-n-1,n+1]$ which satisfy $s\equiv n+1\ \mathrm{mod}\ 2$ and $s\equiv (-1)^n \ \mathrm{mod}\ 3$, and
$c_{n+1,s}$ are entries in the signed Catalan triangle, is indeed the same as $A_n(q)$.
It would be interesting to see what other combinatorics may arise from considering the sum of $F_\lambda(q)$ over $\lambda \vdash n$ with at most $k$ columns for a fixed $k$.
\end{remark}


\section{Jordan canonical forms and $q$-rook placements}\label{sec.rooks}

In light of Corollary~\ref{cor.multiplicity}, we define polynomials $G_\lambda(q) \in \mathbb{Z}[q]$ by
\begin{equation}\label{eqn.Glambda}
F_\lambda(q) = (q-1)^{n-\ell(\lambda)}G_\lambda(q).
\end{equation}
In fact, we can deduce from Corollary~\ref{cor.multiplicity} that $G_\lambda(q)\in \mathbb{N}[q]$.
In this section, we explore the connection between the nonnegative coefficients of $G_\lambda(q)$ and rook placements.

\subsection{Background on rook polynomials}
A {\em board} $B$ is a subset of an $n$ by $n$ grid of squares.  In this paper, we follow Haglund~\cite{haglund98} and Solomon~\cite{solomon90}, and index the squares using the convention for the entries of a matrix.  A {\em Ferrers board} is a board $B$ where if a square $s\in B$, then every square lying north and/or east of $s$ is also in $B$.  Our Ferrers boards have squares justified upwards and to the right.  
Let $B_n$ denote the staircase-shaped board with $n$ columns of sizes $0,1,\ldots, n-1$. Let $\mathrm{area}(B)$ be the number of squares in $B$, so that in particular, $\mathrm{area}(B_n) = {n \choose 2}$.

A placement of $k$ rooks on a board $B$ is {\em non-attacking} if there is at most one rook in each row and each column of $B$.
Let $\mathcal{C}(B,k)$ denote the set of non-attacking placements of $k$ rooks on $B$.  All rook placements considered in this article are non-attacking, so from this point forward, we drop the qualifier.
For a placement $C\in \mathcal{C}(B,k)$, let $\mathrm{ne}(C)$ be the number of squares in $B$ lying directly north or directly east of a rook. The {\em inversion} of the placement is the number
\begin{equation}\label{eq.inv}
\mathrm{inv}(C) = \mathrm{area}(B) - k - \mathrm{ne}(C).
\end{equation}
As noted in~\cite{garsiaRemmel86}, the statistic $\mathrm{inv}(C)$ is a generalization of the number of inversions of a permutation, since permutations can be identified with rook placements on a square-shaped board.


For $i=1,\ldots,n$, the weight of the $i$th column $C_i$ of $C$ is
\begin{equation}
C_i(q) = (q-1)^{\#\mathrm{rooks\ in\ }C_i} q^{\mathrm{ne}(C_i)},
\end{equation}
and the {\em weight} of $C$ is defined by $F_C(q) = \prod_{i=1}^n C_i(q)$.  Alternatively, if $C\in \mathcal{C}(B,k)$, then $F_C(q) = (q-1)^kq^{\mathrm{ne}(C)}$.  


\begin{example}\label{eg.three}
We use $\times$ to mark a rook and use {\tiny$\,\bullet\,$} to mark squares lying directly north or directly east of a rook (these squares shall be referred to as the north-east squares of the placement). 
The following illustration is a placement of four rooks on the staircase-shaped board $B_7$.
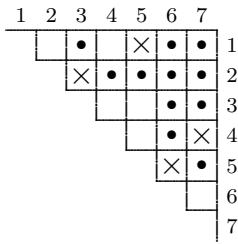
\begin{displaymath}
\beginpicture
\setcoordinatesystem units <.4cm,.4cm> 
\setplotarea x from -7 to 0, y from -7 to 0 
\plot -7 0 0 0 / \plot -6 -1 0 -1 / \plot -5 -2 0 -2 / \plot -4 -3 0 -3 /
\plot -3 -4 0 -4 / \plot -2 -5 0 -5 / \plot -1 -6 0 -6 / 
\plot -6 0 -6 -1 /
\plot -5 0 -5 -2 / \plot -4 0 -4 -3 / \plot -3 0 -3 -4 / \plot -2 0 -2 -5 /
\plot -1 0 -1 -6 / \plot 0 0 0 -7 /
\put{$\times$} at -4.5 -1.5 \put{$\times$} at -2.5 -0.5 
\put{$\times$} at -1.5 -4.5 \put{$\times$} at -0.5 -3.5
\put{\tiny$\bullet$} at -4.5 -0.5 
\put{\tiny$\bullet$} at -3.5 -1.5 \put{\tiny$\bullet$} at -2.5 -1.5
\put{\tiny$\bullet$} at -1.5 -3.5 \put{\tiny$\bullet$} at -1.5 -2.5      
\put{\tiny$\bullet$} at -1.5 -1.5 \put{\tiny$\bullet$} at -1.5 -0.5  
\put{\tiny$\bullet$} at -0.5 -0.5 \put{\tiny$\bullet$} at -0.5 -1.5
\put{\tiny$\bullet$} at -0.5 -2.5 \put{\tiny$\bullet$} at -0.5 -4.5      
\put{\tiny$1$} at -6.5 0.5
\put{\tiny$2$} at -5.5 0.5
\put{\tiny$3$} at -4.5 0.5
\put{\tiny$4$} at -3.5 0.5
\put{\tiny$5$} at -2.5 0.5
\put{\tiny$6$} at -1.5 0.5 
\put{\tiny$7$} at -0.5 0.5
\put{\tiny$1$} at 0.5 -0.5
\put{\tiny$2$} at 0.5 -1.5
\put{\tiny$3$} at 0.5 -2.5
\put{\tiny$4$} at 0.5 -3.5
\put{\tiny$5$} at 0.5 -4.5
\put{\tiny$6$} at 0.5 -5.5
\put{\tiny$7$} at 0.5 -6.5
\endpicture
\end{displaymath}
This rook placement has $\mathrm{ne}(C) = 11$, $\mathrm{inv}(C) = 6$, and weight $F_C(q) = (q-1)^4q^{11}$.
\end{example}

For $k\geq0$, the {\em $q$-rook polynomial} of a Ferrers board $B$ is defined by Garsia and Remmel~\cite[I.4]{garsiaRemmel86} as
\begin{equation}\label{eq.qrook}
R_{B,k}(q) = \sum_{C\in \mathcal{C}(B,k)} q^{\mathrm{inv}(C)}.
\end{equation}

The following result explains the role of rook polynomials in the enumeration of matrices of given rank. 
The {\em support} of a matrix $X$ is $\{(i,j) \mid x_{ij} \neq0\}$.  Given a Ferrers board $B$ with $n$ columns, we may identify the squares in $B$ with the entries in an $n$ by $n$ matrix.
\begin{theorem}[Haglund] \label{thm.PkB}
If $B$ is a Ferrers board, then the number $P_{B,k}(q)$ of $n$ by $n$ matrices of rank $k$ with support contained in $B$ is
\begin{displaymath}
P_{B,k}(q)= (q-1)^{k} q^{\mathrm{area}(B)-k}R_{B,k}(q^{-1}).
\end{displaymath}
\end{theorem}
Looking ahead, it will be convenient to consider Theorem~\ref{thm.PkB} in the following equivalent form:
\begin{equation}\label{eqn.PkB2}
P_{B,k}(q) 
= \sum_{C\in\mathcal{C}(B,k)} (q-1)^{k} q^{\mathrm{ne}(C)}
= \sum_{C\in\mathcal{C}(B,k)} F_C(q). 
\end{equation}

\begin{example}
We list the seven rook placements on $B_4$ with two rooks, along with their weights.
\begin{equation}\label{eq.line1}
\beginpicture 
\setcoordinatesystem units <2.5cm,1cm>
\setplotarea x from 0 to 0, y from 0 to 1.2
\put{\beginpicture
\setcoordinatesystem units <.4cm,.4cm> 
\setplotarea x from -5 to 0, y from -4 to 0
\plot -4 0 0 0 / \plot -3 -1 0 -1 / \plot -2 -2 0 -2 / \plot -1 -3 0 -3 /
\plot -3 0 -3 -1 / \plot -2 0 -2 -2 / \plot -1 0 -1 -3 / \plot 0 0 0 -4 /
\put{$\times$} at -2.5 -0.5 
\put{$\times$} at -1.5 -1.5 
\put{\tiny$\bullet$} at -1.5 -0.5 \put{\tiny$\bullet$} at -0.5 -0.5
\put{\tiny$\bullet$} at -0.5 -1.5
\put{$(q-1)^2q^3$} at -3 -4
\endpicture} at 0 0
\put{\beginpicture
\setcoordinatesystem units <.4cm,.4cm> 
\setplotarea x from -5 to 0, y from -4 to 0
\plot -4 0 0 0 / \plot -3 -1 0 -1 / \plot -2 -2 0 -2 / \plot -1 -3 0 -3 /
\plot -3 0 -3 -1 / \plot -2 0 -2 -2 / \plot -1 0 -1 -3 / \plot 0 0 0 -4 /
\put{$\times$} at -2.5 -0.5 
\put{$\times$} at -0.5 -2.5 
\put{\tiny$\bullet$} at -1.5 -0.5 
\put{\tiny$\bullet$} at -0.5 -0.5
\put{\tiny$\bullet$} at -0.5 -1.5
\put{$(q-1)^2q^3$} at -3 -4
\endpicture} at 1 0
\put{\beginpicture
\setcoordinatesystem units <.4cm,.4cm> 
\setplotarea x from -4 to 0, y from -4 to 0
\plot -4 0 0 0 / \plot -3 -1 0 -1 / \plot -2 -2 0 -2 / \plot -1 -3 0 -3 /
\plot -3 0 -3 -1 / \plot -2 0 -2 -2 / \plot -1 0 -1 -3 / \plot 0 0 0 -4 /
\put{$\times$} at -1.5 -1.5 
\put{$\times$} at -0.5 -2.5 
\put{\tiny$\bullet$} at -1.5 -0.5 
\put{\tiny$\bullet$} at -0.5 -0.5
\put{\tiny$\bullet$} at -0.5 -1.5
\put{$(q-1)^2q^3$} at -3 -4
\endpicture} at 2 0
\put{\beginpicture
\setcoordinatesystem units <.4cm,.4cm> 
\setplotarea x from -4 to 0, y from -4 to 0
\plot -4 0 0 0 / \plot -3 -1 0 -1 / \plot -2 -2 0 -2 / \plot -1 -3 0 -3 /
\plot -3 0 -3 -1 / \plot -2 0 -2 -2 / \plot -1 0 -1 -3 / \plot 0 0 0 -4 /
\put{$\times$} at -1.5 -0.5 
\put{$\times$} at -0.5 -2.5 
\put{\tiny$\bullet$} at -0.5 -0.5
\put{\tiny$\bullet$} at -0.5 -1.5
\put{$(q-1)^2q^2$} at -3 -4
\endpicture} at 3 0
\endpicture
\end{equation}

\begin{equation}\label{eq.line2}
\beginpicture
\setcoordinatesystem units <2.5cm,1cm> 
\setplotarea x from 0 to 0, y from -1 to 0
%
%
\put{\beginpicture
\setcoordinatesystem units <.4cm,.4cm> 
\setplotarea x from -4 to 0, y from -4 to 0
\plot -4 0 0 0 / \plot -3 -1 0 -1 / \plot -2 -2 0 -2 / \plot -1 -3 0 -3 /
\plot -3 0 -3 -1 / \plot -2 0 -2 -2 / \plot -1 0 -1 -3 / \plot 0 0 0 -4 /
\put{$\times$} at -2.5 -0.5 
\put{$\times$} at -0.5 -1.5 
\put{\tiny$\bullet$} at -1.5 -0.5 
\put{\tiny$\bullet$} at -0.5 -0.5
\put{$(q-1)^2q^2$} at -3 -4
\endpicture} at 0 0
\put{\beginpicture
\setcoordinatesystem units <.4cm,.4cm> 
\setplotarea x from -4 to 0, y from -4 to 0
\plot -4 0 0 0 / \plot -3 -1 0 -1 / \plot -2 -2 0 -2 / \plot -1 -3 0 -3 /
\plot -3 0 -3 -1 / \plot -2 0 -2 -2 / \plot -1 0 -1 -3 / \plot 0 0 0 -4 /
\put{$\times$} at -0.5 -0.5 
\put{$\times$} at -1.5 -1.5 
\put{\tiny$\bullet$} at -1.5 -0.5
\put{\tiny$\bullet$} at -0.5 -1.5
\put{$(q-1)^2q^2$} at -3 -4
\endpicture} at 1 0
\put{\beginpicture
\setcoordinatesystem units <.4cm,.4cm> 
\setplotarea x from -4 to 0, y from -4 to 0
\plot -4 0 0 0 / \plot -3 -1 0 -1 / \plot -2 -2 0 -2 / \plot -1 -3 0 -3 /
\plot -3 0 -3 -1 / \plot -2 0 -2 -2 / \plot -1 0 -1 -3 / \plot 0 0 0 -4 /
\put{$\times$} at -1.5 -0.5 
\put{$\times$} at -0.5 -1.5 
\put{\tiny$\bullet$} at -0.5 -0.5
\put{$(q-1)^2q$} at -3 -4
\endpicture} at 2 0
\endpicture
\end{equation}
Thus 
$P_{B_4,2}(q) = (q-1)^2(3q^3+3q^2+q)$.
\end{example}

\subsection{Rook placements and Jordan forms}

The purpose of this section is to generalize Haglund's formula~\eqref{eqn.PkB2} to a formula for $F_\lambda(q)$ (Corollary~\ref{cor.main}) as a sum over a set of rook placements.  We achieve this by defining a multigraph $\mathcal{Z}$ that is related to $\mathcal{Y}$, and show that paths in $\mathcal{Z}$ are equivalent to rook placements.

The multigraph $\mathcal{Z}$ is constructed from $\mathcal{Y}$ by replacing each edge of $\mathcal{Y}$ by one or more edges as follows.
If there is an edge from $\mu$ to $\lambda$ in $\mathcal{Y}$ of weight 
$q^{|\mu|-\mu_{j-1}'}\left(q^{\mu_{j-1}'-\mu_j'}-1 \right)$, then 
this edge is replaced by $\mu_{j-1}'-\mu_j'$ edges from $\mu$ to $\lambda$ with weights
\begin{equation}\label{eqn.deg}
(q-1)q^{|\mu|-\mu_j'-1}, \ldots, (q-1)q^{|\mu|-\mu_{j-1}'}
\end{equation} 
in $\mathcal{Z}$.  All other edges remain as before.  See Figure~\ref{fig.Z}.



\begin{figure}
\Yboxdim8pt
$$\xymatrixcolsep{4pc}\xymatrixrowsep{3pc}
\xymatrix{
& {\yng(1,1,1,1)}
	\\
& {\yng(1,1,1)} \ar@/^1pc/[r]^{(q-1)q^2} \ar[r]|-{(q-1)q} 
				\ar@/_1pc/[r]_{q-1}		\ar[u]_-1
	& {\yng(2,1,1)} 
	& {\yng(2,2)} 
	\\
& {\yng(1,1)} \ar@/^/[r]^{(q-1)q} \ar@/_/[r]_{q-1}			\ar[u]_-1
	& {\yng(2,1)} \ar[r]^{(q-1)q^2}  \ar[u]_-q 		\ar[ur]|{(q-1)q}
	& {\yng(3,1)} 
	\\
& {\yng(1)} \ar[r]^{q-1} 			\ar[u]_-1
	& {\yng(2)} \ar[r]^{(q-1)q} 	\ar[u]_-q
	& {\yng(3)} \ar[r]^{(q-1)q^2} 	\ar[u]_-{q^2}
	& {\yng(4)} 
	\\
\emptyset \ar[ur]_1 
}
$$
\caption{The multigraph $\mathcal{Z}$, up to $n=4$.}
\label{fig.Z}
\end{figure}

Let $\mathcal{P}_{\mathcal{Z}}(\lambda)$ denote the set of paths in the graph $\mathcal{Z}$ from the empty partition $\emptyset$ to $\lambda$.
For a path $P = (\emptyset = \pi^{(0)}, \pi^{(1)}, \ldots, \pi^{(n)}=\lambda)$ in $\mathcal{P}_\mathcal{Z}(\lambda)$, let $\epsilon_i(q)$
denote the weight of the $i$th edge, for $i=1,\ldots,n$.  Naturally, we define the weight of the path by $F_P(q) = \prod_{i=1}^n \epsilon_i(q)$,
so that
\begin{equation}\label{eq.formula}
F_\lambda(q) = \sum_{P\in \mathcal{P}_{\mathcal{Z}}(\lambda)} F_P(q).
\end{equation}

\begin{lemma}\label{lem.extension}
Let $\mu\vdash n-1$ be a partition with $\ell(\mu) = \ell$ parts.  Then there are $\ell+1$ edges leaving $\mu$ in the graph $\mathcal{Z}$, with weights
$$(q-1)q^{|\mu|-1}, (q-1)q^{|\mu|-2}, \ldots,  (q-1)q^{|\mu|-\ell}, \hbox{ and } q^{|\mu|-\ell}.$$
\end{lemma}
\proof
If a partition $\lambda \vdash n$ is obtained by adding a box to the first column of $\mu$, then 
there is a unique edge from $\mu$ to $\lambda$ in $\mathcal{Z}$ with weight $q^{|\mu|-\ell}$.
Otherwise, if we consider the set of all partitions which can be obtained from $\mu$ by adding a box anywhere except in the first column, then there are a total of
$$\sum_{j\geq2} \left( \mu_{j-1}'-\mu_j' \right) = \ell $$
edges from $\mu$ to some partition of $n$.  Moreover, by Equation~\eqref{eqn.deg}, these $\ell$ weights are $(q-1) q^{|\mu|-i}$ for $i=1,\ldots, \ell$.
\qed

A sequence of nonnegative integers is $\mathcal{P}_\mathcal{Z}$-admissible if it is the degree sequence of a path $P=(\emptyset, \pi^{(1)},\ldots, \pi^{(n)})$ in $\mathcal{Z}$.  That is, $(d_1,\ldots, d_n) = (\deg\epsilon_1(q),\ldots, \deg\epsilon_n(q))$.

\begin{corollary}\label{cor.degsequence}
A $\mathcal{P}_\mathcal{Z}$-admissible sequence determines a unique path in $\mathcal{Z}$.
\end{corollary}
\proof Induct on $n$.  When $n=1$, the only path is the from $\emptyset$ to $(1)$, and it has degree sequence $(0)$.

Given a $\mathcal{P}_\mathcal{Z}$-admissible sequence $(d_1,\ldots, d_n)$, the subsequence $(d_1,\ldots, d_{n-1})$ determines a unique path $P'= (\emptyset, \pi^{(1)},\ldots, \pi^{(n-1)})$.  Suppose $\mu = \pi^{(n-1)}$ has $\ell$ parts.  Then $|\mu|-\ell+1 \leq d_n \leq |\mu|$, and by Lemma~\ref{lem.extension}, there is a unique edge leaving $\mu$ with degree $d_n$.
\qed

\subsection{The construction of $\Phi$}
Let $\mathcal{P}_\mathcal{Z}(n,n-k)$ denote the set of paths in $\mathcal{Z}$ from $\emptyset$ to a partition of $n$ with $n-k$ parts.  In this section, we define a weight-preserving bijection $\Phi: \mathcal{C}(B_n,k)\rightarrow \mathcal{P}_\mathcal{Z}(n,n-k)$.

\begin{prop}\label{prop.main}
Let $n\geq1$ and $k=0,\ldots,n-1$. 
Let $C \in \mathcal{C}(B_n,k)$ be a rook placement with columns $C_1,\ldots, C_n$.  There exists a unique path $P\in \mathcal{P}_\mathcal{Z}(n,n-k)$ with edge weights $(\epsilon_1(q),\ldots, \epsilon_n(q)) = (C_1(q),\ldots, C_n(q))$.
\end{prop}
\proof
Proceed by induction on $n+k$. When $n=1$ and $k=0$, there is a unique rook placement on the empty board $B_1$ with no rooks having weight one, corresponding to the unique path $P = (\emptyset, (1))$ in $\mathcal{Z}$ with the same weight.

Assume the result holds for all rook placements in $\mathcal{C}(B_{n-1}, k)$ and $\mathcal{C}(B_{n-1}, k-1)$.  Given a rook placement $C \in \mathcal{C}(B_n, k)$, let $C'$ be the sub-placement consisting of the first $n-1$ columns of $C$.  By induction, the sequence $(C_1(q),\ldots, C_{n-1}(q))$ determines a unique path $(\emptyset, \pi^{(1)}, \ldots, \pi^{(n-1)})$ in $\mathcal{Z}$ such that $\epsilon_i(q)=C_i(q)$ for $i=1,\ldots, n-1$.

There are now two cases two consider.  The first case is if $C'\in \mathcal{C}(B_{n-1}, k)$, so that $\ell(\pi^{(n-1)})=n-k-1$.  
There are $k$ rooks in $C'$, so the $n$th column of $C$ does not contain any rooks, and $C_n(q) = q^k$.
By Lemma~\ref{lem.extension}, there exists a unique edge in the graph $\mathcal{Z}$ originating at $\pi^{(n-1)}$ with weight $q^k$.  
Thus $C$ corresponds to the path $P= (\emptyset,\pi^{(1)}, \ldots,  \pi^{(n-1)}, \pi^{(n)})$
where $\pi^{(n)}$ is obtained from $\pi^{(n-1)}$ by adding a box to the first column, and $\epsilon_n(q) = q^{k}$.  Moreover, $\ell(\pi^{(n)})= n-k$.

The second case is if $C'\in \mathcal{C}(B_{n-1}, k-1)$, so that 
$\ell(\pi^{(n-1)})=n-k$.  There must be $k-1$ `northeast' squares in the $n$th column of $C$, and there are $n-k$ remaining squares in that column where a rook may be placed.  
Label these available squares $a_0, a_1,\ldots, a_{n-k-1}$ from the top to the bottom.  Observe that $C_n(q) = (q-1)q^{k-1+i}$ if a rook is placed in the square $a_i$, for $0\leq i\leq n-k-1$.
Again by Lemma~\ref{lem.extension}, there exists $n-k$ edges in the graph $\mathcal{Z}$ originating at $\pi^{(n-1)}$ with the weights $(q-1)q^{h}$ for $k-1 \leq h \leq n-2$.  Thus if the $k$th rook of $C$ is placed in the square $a_i$, then $C$ corresponds to the path 
$P= (\emptyset,\pi^{(1)}, \ldots,  \pi^{(n-1)}, \pi^{(n)})$ with $\epsilon_n(q) = (q-1)q^{k-1+i}$, and $\ell(\pi^{(n)}) = n-k$.
\qed

Given $C\in \mathcal{C}(B_n,k)$, let $\Phi(C)$ be the path in $\mathcal{P}_\mathcal{Z}(n,n-k)$ with edge weights $(\epsilon_1(q),\ldots, \epsilon_n(q)) = (C_1(q),\ldots, C_n(q))$.

\begin{theorem}\label{thm.Phi} 
The map $\Phi: \mathcal{C}(B_n,k) \rightarrow \mathcal{P}_\mathcal{Z}(n,n-k)$ is a weight-preserving bijection.
\end{theorem}  
\proof
Proposition~\ref{prop.main} shows that the map $\Phi$ is an injective weight-preserving map, since each column of the rook placement determines each edge of the path $\Phi(C)$:
$$F_C(q) = \prod_{i=1}^n C_i(q) = \prod_{i=1}^n \epsilon_i(q) = F_{\Phi(C)}(q).$$ 
In fact, the proof of the Proposition also shows that $\Phi$ is surjective because the number of possible ways to add a column to an existing rook placement is equal to the number of possible ways to extend a path in $\mathcal{Z}$ by one edge. Therefore, $\Phi$ is a weight-preserving bijection.
\qed


A sequence of nonnegative integers is {\em $\mathcal{C}$-admissible} if it is the degree sequence of a rook placement.  That is, $(d_1,\ldots, d_n)=(\deg C_1(q)), \ldots, \deg C_n(q))$ for a $C \in \mathcal{C}(B_n,k)$.
The next Corollary follows easily from Theorem~\ref{thm.Phi}.

\begin{corollary} \label{cor.cadmissible}
A $\mathcal{C}$-admissible sequence determines a unique rook placement. \qed
\end{corollary}
%

It follows from Theorem~\ref{thm.Phi} that we may associate a partition type to each rook placement on $B_n$.  
The {\em partition type} of a rook placement $C$
is the partition at the endpoint of the path $\Phi(C)$ in $\mathcal{Z}$.
Let $\mathcal{C}(\lambda)=\Phi^{-1}\left(P_{\mathcal{Z}}(\lambda)\right)$ denote the set of rook placements of partition type $\lambda$.

\begin{corollary} \label{cor.main}
Let $\lambda \vdash n$ be a partition with $\ell(\lambda) = n-k$ parts.  Then
$$F_\lambda(q)
= \sum_{C\in \mathcal{C}(\lambda)} F_C(q) 
= (q-1)^{n-\ell(\lambda)} \sum_{C\in \mathcal{C}(\lambda)}  q^{\mathrm{ne}(C)}. $$
\end{corollary}
\proof
The result follows from Equation~\ref{eq.formula} and the bijection $\Phi$.
\qed

\begin{remark}\label{rem.Glambda}
The polynomial $G_\lambda(q)\in \mathbb{N}[q]$ defined in Equation~\eqref{eqn.Glambda} is simply a sum over the rook placements of type $\lambda$ involving the north-east statistic.
\end{remark}

\section{A connection with set partitions}

The results of the previous section naturally leads to a decomposition of $F_T(q)$, indexed by some tableau $T$, into a sum of polynomials indexed by set partitions, which we explain below. 

A {\em set partition} is a set $S=\{s_1, \ldots, s_k\}$ of nonempty disjoint subsets of $[n]$ such that $\bigcup_{i=1}^k s_i= [n]$.  The $s_i$'s are the {\em blocks} of $\sigma$.  Let $\ell(S)$ denote the number of blocks of $S$, and let $\mathcal{S}(n,n-k)$ denote the set of set partitions of $[n]$ with $n-k$ blocks.
We adopt the convention of listing the blocks in order so that 
\begin{equation}\label{eqn.convention}
|s_1|\geq |s_2|\geq \cdots \geq |s_k|, \hbox{ and } 
\min s_i < \min s_{i+1} \hbox{ if } |s_i|=|s_{i+1}|.
\end{equation}  
This allows us to represent a set partition with a diagram similar to that of a standard Young tableau; the $i$th row of the diagram consists of the elements in the block $s_i$ listed in increasing order, but there are no restrictions on the entries in each column of the diagram.  
A set partition $S=(s_1,\ldots, s_m)$ has {\em partition type} $\lambda$ if $\lambda = (|s_1|, \ldots, |s_m|)$.

For $i=1,\ldots,n$, let $S^{(i)}$ denote the sub-diagram of $S$ consisting of the boxes containing $1,\ldots, i$, with rows ordered according to the convention set forth in Equation~\eqref{eqn.convention}. 
If the box containing $i$ is not in the first column of the diagram, let $u$ be the least element in the same row as $i$ in $S^{(i)}$, and suppose $u$ is in the $r$th row of $S^{(i-1)}$ for some $1\leq r\leq \ell(S^{(i-1)})$.
The weight arising from the $i$th box is
\begin{equation}\label{eqn.Sweight}
S^{(i)}(q) = \begin{cases}
q^{i-1-\ell(S^{(i-1)})}, 
	& \hbox{if the $i$th box is in the first column,}\\
(q-1)q^{i-1-r}, 
 	& \hbox{if the $i$th box is in the $j$th column, $j\geq2$.}
\end{cases}
\end{equation}
We define the {\em weight} of $S$ as $F_S(q) = \prod_{i=1}^n S^{(i)}(q)$.  

A sequence of nonnegative integers is {\em $\mathcal{S}$-admissible} if it is the degree sequence of a set partition.  That is, $(d_1,\ldots, d_n)=(\deg S^{(1)}(q)), \ldots, \deg S^{(n)}(q))$ for a $S \in \mathcal{S}(n)$.

\begin{lemma}\label{lem.sadmissible} An $\mathcal{S}$-admissible sequence determines a unique set partition.
\end{lemma}
\proof Induct on $n$. When $n=1$, the only set partition is $\{\{1\}\}$, and its degree sequence is $(0)$.

Given an $\mathcal{S}$-admissible sequence $(d_1,\ldots, d_n)$, the subsequence $(d_1,\ldots, d_{n-1})$ determines a unique set partition $S^{(n-1)}= (S^{(n-1)}_1,\ldots, S^{(n-1)}_m)$.
By Equation~\eqref{eqn.Sweight},
$n-1-m \leq d_n \leq n-1$, and each of the $m+1$ choices for $d_n$ determines the block of $S^{(n-1)}$ into which $n$ should be inserted. \qed

We have already constructed a weight-preserving bijection $\Phi$ between rook placements and paths in $\mathcal{Z}$.  We now construct a weight-preserving bijection $\Psi$ between rook placements and set partitions, effectively showing that paths in $\mathcal{Z}$ are equivalent to set partitions, so that  
$F_Z(q) = F_C(q) = F_S(q)$ 
if 
$Z \longleftrightarrow C \longleftrightarrow S$
for $Z\in\mathcal{P}_\mathcal{Z}(n,n-k)$, $C\in \mathcal{C}(n,k)$, and $S\in \mathcal{S}(n,n-k)$.

\begin{remark}
There is a classically known bijection (see~\cite{stanley99}) between the set of rook placements on the staircase board $B_n$ with $k$ rooks and the set of set partitions of $[n] = \{1,\ldots, n\}$ with $n-k$ blocks: the placement $C$ corresponds to the set partition where the integers $i$ and $j$ are in the same block if and only if there is a rook in the square $(i,j)\in C$.  This bijection is different from the one described in Theorem~\ref{thm.Psi}.  For example, the classical bijection associates the rook placement
$$
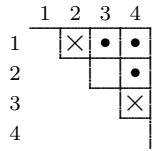
\beginpicture
\setcoordinatesystem units <.4cm,.4cm> 
\setplotarea x from -4 to 0, y from -4 to 0
\plot -4 0 0 0 / \plot -3 -1 0 -1 / \plot -2 -2 0 -2 / \plot -1 -3 0 -3 /
\plot -3 0 -3 -1 / \plot -2 0 -2 -2 / \plot -1 0 -1 -3 / \plot 0 0 0 -4 /
\put{$\times$} at -2.5 -0.5 
\put{$\times$} at -0.5 -2.5 
\put{\tiny$\bullet$} at -1.5 -0.5 
\put{\tiny$\bullet$} at -0.5 -0.5
\put{\tiny$\bullet$} at -0.5 -1.5
\multiput{\tiny$1$} at -3.5 .5 -4.5 -.5 /
\multiput{\tiny$2$} at -2.5 .5 -4.5 -1.5 /
\multiput{\tiny$3$} at -1.5 .5 -4.5 -2.5 /
\multiput{\tiny$4$} at -0.5 .5 -4.5 -3.5 /
\endpicture
$$
to the set partition $(\{1,2\},\{3,4\})$ and so 
has partition type $(2,2)$, but as we shall see below, this placement is associated to the set partition $(\{1,2,4\}, \{3\})$ under the bijection in Theorem~\ref{thm.Psi} and has partition type $(3,1)$.
\end{remark}

\subsection{The construction of $\Psi$}

Let $C\in \mathcal{C}(B_n,k)$ be a rook placement. The main idea is that the degree of $C_i(q)$ arising from the $i$th column of $C$ determines the block of the set partition in which we place $i$.  In the construction of the set partition $\Psi(C)$, we will create a sequence of intermediate set partitions $S^{(i)}$ of $[i]$ for $i=1,\ldots, n$.

The initial case is always $\deg(C_1(q))=  \deg(1)= 0$, so $S^{(1)} = \{\{1\}\}$. Assume that $S^{(i-1)} = \{S^{(i-1)}_1 , \ldots , S^{(i-1)}_m\}$ is the set partition which corresponds to the first $i-1$ columns of $C$, so that $m=\ell(S^{(i-1)})$.   
Observe that there are $m+1$ possible blocks in which to insert $i$ to obtain $S^{(i)}$.
By Corollary~\ref{cor.degsequence}, we know that 
$$i-1-\ell(S^{(i-1)}) \leq \deg(C_i(q)) \leq i-1,$$
so we construct $S^{(i)}$ by placing $i$ in the $j$th block of $S^{(i-1)}$, where $j= i-\deg(C_i(q))$, and then rearranging the blocks to fit the convention in Equation~\eqref{eqn.convention} if necessary.

\begin{theorem}\label{thm.Psi}
The map $\Psi: \mathcal{C}(n,k) \rightarrow \mathcal{S}(n,n-k)$ is a weight-preserving bijection.
\end{theorem}
\proof  Let $S=\Psi(C)$.  The map $\Psi$ is weight-preserving, as $C_i(q) = S^{(i)}(q)$ by construction, for each $i=1,\ldots, n$.
Now, since the degrees $\deg C_i(q)= \deg S^{(i)}(q)$, and by Corollary~\ref{cor.cadmissible} and Lemma~\ref{lem.sadmissible} the sequences of degrees completely determine $C$ and $S$ respectively, then $\Psi$ is injective.  
Finally, we note that $|\mathcal{C}(n,k)|=|\mathcal{S}(n,n-k)|$, so $\Psi$ is a bijection. \qed


\begin{corollary}
Let $\mathcal{S}(\lambda)$ denote the set of all set partitions of partition type $\lambda$.  Then
$$F_\lambda(q) = \sum_{S\in \mathcal{S}(\lambda)} F_S(q). $$
\qed
\end{corollary}

\begin{example} Let $C$ be
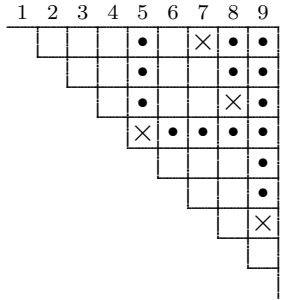
\begin{displaymath}
\beginpicture
\setcoordinatesystem units <.4cm,.4cm> 
\setplotarea x from -9 to 0, y from -9 to 1
\plot -9 0 0 0 / \plot -8 -1 0 -1 / \plot -7 -2 0 -2 / \plot -6 -3 0 -3 / 
\plot -5 -4 0 -4 / \plot -4 -5 0 -5 / \plot -3 -6 0 -6 / \plot -2 -7 0 -7 / 
\plot -1 -8 0 -8 /

\plot -8 0 -8 -1 / \plot -7 0 -7 -2 / \plot -6 0 -6 -3 / \plot -5 0 -5 -4 /
\plot -4 0 -4 -5 / \plot -3 0 -3 -6 / \plot -2 0 -2 -7 / \plot -1 0 -1 -8 / \plot 0 0 0 -9 /

\put{$\times$} at -4.5 -3.5 
\put{$\times$} at -2.5 -0.5 
\put{$\times$} at -1.5 -2.5 
\put{$\times$} at -0.5 -6.5 

\multiput{\tiny$\bullet$} at -4.5 -0.5 -4.5 -1.5 -4.5 -2.5 /
\put{\tiny$\bullet$} at -3.5 -3.5 
\multiput{\tiny$\bullet$} at -2.5 -3.5 /
\multiput{\tiny$\bullet$} at -1.5 -0.5 -1.5 -1.5 -1.5 -3.5 /
\multiput{\tiny$\bullet$} at -0.5 -0.5 -0.5 -1.5 -0.5 -2.5 -0.5 -3.5 -0.5 -4.5 -0.5 -5.5 /

\put{\tiny$1$} at -8.5 0.5 \put{\tiny$2$} at -7.5 0.5 
\put{\tiny$3$} at -6.5 0.5 \put{\tiny$4$} at -5.5 0.5 
\put{\tiny$5$} at -4.5 0.5 \put{\tiny$6$} at -3.5 0.5
\put{\tiny$7$} at -2.5 0.5 \put{\tiny$8$} at -1.5 0.5
\put{\tiny$9$} at -0.5 0.5 
\endpicture
\end{displaymath}
The associated sequence of set partition diagrams associated to $C$ is
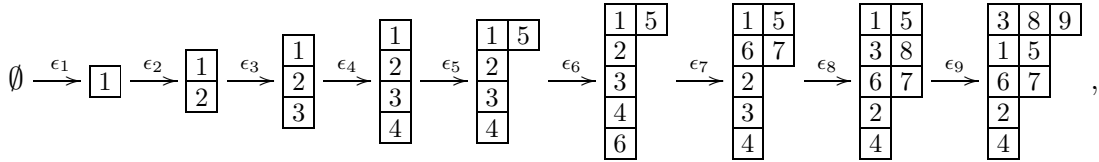
\Yboxdim12pt
\begin{displaymath}
\xymatrix@C=1.5em{
{\footnotesize\emptyset} \ar[r]^{\epsilon_1} &
{\footnotesize\young(1)} \ar[r]^{\epsilon_2} &
{\footnotesize\young(1,2)} \ar[r]^{\epsilon_3} &
{\footnotesize\young(1,2,3)} \ar[r]^{\epsilon_4} &
{\footnotesize\young(1,2,3,4)} \ar[r]^{\epsilon_5} &
{\footnotesize\young(15,2,3,4)} \ar[r]^{\epsilon_6} &
{\footnotesize\young(15,2,3,4,6)} \ar[r]^{\epsilon_7} &
{\footnotesize\young(15,67,2,3,4)} \ar[r]^{\epsilon_8} &
{\footnotesize\young(15,38,67,2,4)} \ar[r]^-{\epsilon_9} &
{\footnotesize\young(389,15,67,2,4)} 
},
\end{displaymath}
so the set partition associated to the rook placement $C$ is
$$S=\Psi(C) = (\{3,8,9\}, \{1,5\}, \{6,7\}, \{2\}, \{4\}).$$
\end{example}

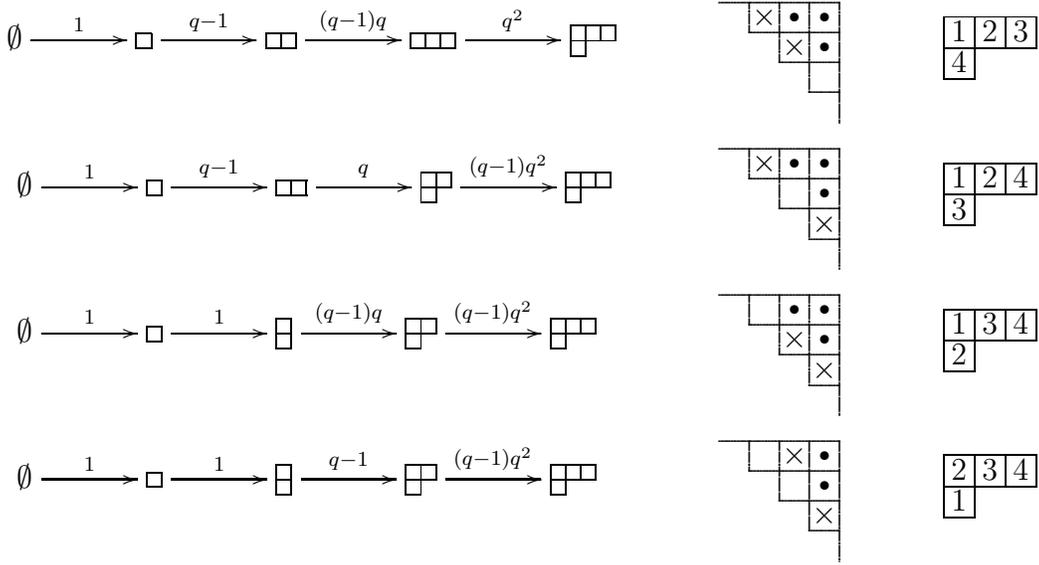
\begin{figure}\label{eg.mother}
\Yboxdim6pt
\begin{tabular}{c}
\beginpicture
	\setcoordinatesystem units <.4cm,.4cm> 
	\setplotarea x from 0 to 0, y from -4 to 0
	
	\put{\xymatrixcolsep{3pc}\xymatrix{\emptyset\ar[r]^-1 
		& {\yng(1)}\ar[r]^-{q-1} 
		& {\yng(2)}\ar[r]^-{(q-1)q} 
		& {\yng(3)} \ar[r]^-{q^2} 
		& {\yng(3,1)}}}[l] at -28 -1
	
	\plot -4 0 0 0 / \plot -3 -1 0 -1 / \plot -2 -2 0 -2 / 
	\plot -1 -3 0 -3 /
	\plot -3 0 -3 -1 / \plot -2 0 -2 -2 / \plot -1 0 -1 -3 / 
	\plot 0 0 0 -4 /
	\multiput{$\times$} at -2.5 -0.5 -1.5 -1.5 /
	\multiput{\tiny$\bullet$} at -1.5 -0.5 -0.5 -0.5 -0.5 -1.5 /
	
	\put{$\Yboxdim12pt\young(123,4)$} at 5 -1.5 
	\endpicture 
\\
\beginpicture
\setcoordinatesystem units <.4cm,.4cm> 
\setplotarea x from 0 to 0, y from -4 to 0
	\put{\xymatrixcolsep{3pc}
	\xymatrix{\emptyset \ar[r]^-1 
	& {\yng(1)} \ar[r]^-{q-1} 
	& {\yng(2)} \ar[r]^-{q} 
	& {\yng(2,1)} \ar[r]^-{(q-1)q^2} 
	& {\yng(3,1)}} }[l] at -28 -1
	
	\plot -4 0 0 0 / \plot -3 -1 0 -1 / \plot -2 -2 0 -2 / 
	\plot -1 -3 0 -3 / \plot -3 0 -3 -1 / \plot -2 0 -2 -2 / 
	\plot -1 0 -1 -3 / \plot 0 0 0 -4 /
	\multiput{$\times$} at -2.5 -0.5 -0.5 -2.5 / 
	\multiput{\tiny$\bullet$} at -1.5 -0.5 -0.5 -0.5 -0.5 -1.5 /

	\put{$\Yboxdim12pt\young(124,3)$} at 5 -1.5
	\endpicture
\\
\beginpicture
	\setcoordinatesystem units <.4cm,.4cm> 
	\setplotarea x from 0 to 0, y from -4 to 0
	
	\put{\xymatrixcolsep{3pc}
	\xymatrix{\emptyset \ar[r]^-1 & {\yng(1)} \ar[r]^-1 
		& {\yng(1,1)}\ar[r]^-{(q-1)q} 
		& {\yng(2,1)} \ar[r]^-{(q-1)q^2} 
		& {\yng(3,1)}} }[l] at -28 -1
		
	\plot -4 0 0 0 / \plot -3 -1 0 -1 / \plot -2 -2 0 -2 / 
	\plot -1 -3 0 -3 / \plot -3 0 -3 -1 / \plot -2 0 -2 -2 / 
	\plot -1 0 -1 -3 / \plot 0 0 0 -4 /
	\multiput{$\times$} at -1.5 -1.5 -0.5 -2.5 / 
	\multiput{\tiny$\bullet$} at -1.5 -0.5 -0.5 -0.5 -0.5 -1.5 /

	\put{$\Yboxdim12pt\young(134,2)$} at 5 -1.5
	\endpicture
\\
\beginpicture
	\setcoordinatesystem units <.4cm,.4cm> 
	\setplotarea x from 0 to 0, y from -4 to 0
	\put{\xymatrixcolsep{3pc}
		\xymatrix{\emptyset \ar[r]^-1 & {\yng(1)} \ar[r]^-1 
		& {\yng(1,1)} \ar[r]^-{q-1} 
		& {\yng(2,1)} \ar[r]^-{(q-1)q^2} 
		& {\yng(3,1)}} }[l] at -28 -1

	\plot -4 0 0 0 / \plot -3 -1 0 -1 / \plot -2 -2 0 -2 / 
	\plot -1 -3 0 -3 / \plot -3 0 -3 -1 / \plot -2 0 -2 -2 / 
	\plot -1 0 -1 -3 / \plot 0 0 0 -4 /
	\multiput{$\times$} at -1.5 -0.5 -0.5 -2.5 /
	\multiput{\tiny$\bullet$} at -0.5 -0.5 -0.5 -1.5 /

	\put{$\Yboxdim12pt\young(234,1)$} at 5 -1.5
\endpicture
\\
\end{tabular}

\caption{Paths, rook placements, and set partitions related to the computation of
$F_{(3,1)}(q) = (q-1)^2(3q^3+q^2).$
}
\label{fig.mother}
\end{figure}

\begin{remark}
An intriguing question is to ask for a geometric interpretation of the polynomials $F_C(q)$, indexed by rook placements (or set partitions or paths in $\mathcal{Z}$).  

The problem of determining the number of adjoint $G_n(\mathbb{F}_q)$ orbits on $\mathfrak{g}_n(\mathbb{F}_q)$ remains open.  In the case $q=2$, this number has been computed for $n\leq 16$ by Pak and Soffer~\cite[Appendix B]{paksoffer15}. 
Let $\mathcal{O}_{n}(k)$ denote the orbits of rank $k$ matrices.  When $k=1$, it turns out that the polynomials $F_C(q)$ indexed by rook placements with exactly one rook gives the sizes of the ${n\choose 2}$ orbits in $\mathcal{O}_n(1)$. 
For $2\leq i <j \leq n$, each orbit contains a unique matrix $E_{ij}$ whose $ij$th entry is $1$, and is zero everywhere else. The orbit containing $E_{ij}$ is associated to the rook placement $C(i,j)$ with a single rook in the $ij$th square, and the size of the associated orbit is $F_{C(i,j)}(q) = (q-1) q^{n-1-(j-i)}$.

In particular, the formula in Proposition~\ref{prop.hookshapes} applied to the partition $\lambda= (2,1^{n-2})$ 
gives the generating function 
$$F_{(2,1^{n-2})}(q) = (q-1)\left( (n-1)q^{n-2} + (n-2)q^{n-3} +\cdots + 3q^2+2q+1\right) $$
for rank one orbits of $G_n(\mathbb{F}(q))$ on $\mathfrak{g}_n(\mathbb{F}(q))$.
\end{remark}

\begin{remark}
To close, we mention a related problem which may provide a geometrical interpretation of $F_C(q)$ for every rook placement $C$. 
Let $N$ be an $n\times n$ nilpotent matrix with entries in an algebraically closed field $k$ containing $\mathbb{F}_q$, and suppose $N$ has Jordan type $\lambda\vdash n$.  A complete flag $f=(f_1,\ldots, f_n)$ is a sequence of subspaces in $k^n$ such that $f_1\subset \cdots \subset f_n$ and $\dim f_i=i$ for all $i$. A flag is {\em $N$-stable} if $N(f_i) \subseteq f_i$ for all $i$.  Spaltenstein~\cite{spaltenstein76} showed that the variety $X_\lambda$ of $N$-stable flags is a disjoint union of $f^\lambda$ smooth irreducible subvarieties $X_T$ indexed by the standard Young tableaux of shape $\lambda$.  Moreover, the closures $\overline{X}_T$ are the irreducible components of $X_\lambda$, each of which has dimension $n_\lambda$.  The number of $\mathbb{F}_q$-rational points in $X_\lambda$ is given by Green's polynomials $Q_{(1^n)}^\lambda(q)$.  Evidently,
$$\left(\prod_{i\geq1} [m_i(\lambda)]_q! \right)^{-1} Q_{(1^n)}^\lambda(q) = \left((q-1)^{n-\ell(\lambda)} q^m \right)^{-1} F_\lambda(q),$$
with $m=\min_{C\in \mathcal{C}(\lambda)}\mathrm{ne}(C)$.
Based on some computations for small values of $n$, we expect that $F_C(q)$ plays a role in counting points in certain intersections of the irreducible components $\overline{X}_T$.
\end{remark}

\end{document}